\magnification=1200
%\nopagenumbers

\hsize=11.25cm    
\vsize=18cm       
\parindent=12pt   \parskip=5pt     

\hoffset=.5cm   
\voffset=.8cm   

\pretolerance=500 \tolerance=1000  \brokenpenalty=5000

\catcode`\@=11

\font\eightrm=cmr8         \font\eighti=cmmi8
\font\eightsy=cmsy8        \font\eightbf=cmbx8
\font\eighttt=cmtt8        \font\eightit=cmti8
\font\eightsl=cmsl8        \font\sixrm=cmr6
\font\sixi=cmmi6           \font\sixsy=cmsy6
\font\sixbf=cmbx6

\font\tengoth=eufm10 
\font\eightgoth=eufm8  
\font\sevengoth=eufm7      
\font\sixgoth=eufm6        \font\fivegoth=eufm5

\skewchar\eighti='177 \skewchar\sixi='177
\skewchar\eightsy='60 \skewchar\sixsy='60

\newfam\gothfam           \newfam\bboardfam

\def\tenpoint{
  \textfont0=\tenrm \scriptfont0=\sevenrm \scriptscriptfont0=\fiverm
  \def\rm{\fam\z@\tenrm}
  \textfont1=\teni  \scriptfont1=\seveni  \scriptscriptfont1=\fivei
  \def\oldstyle{\fam\@ne\teni}\let\old=\oldstyle
  \textfont2=\tensy \scriptfont2=\sevensy \scriptscriptfont2=\fivesy
  \textfont\gothfam=\tengoth \scriptfont\gothfam=\sevengoth
  \scriptscriptfont\gothfam=\fivegoth
  \def\goth{\fam\gothfam\tengoth}
  
  \textfont\itfam=\tenit
  \def\it{\fam\itfam\tenit}
  \textfont\slfam=\tensl
  \def\sl{\fam\slfam\tensl}
  \textfont\bffam=\tenbf \scriptfont\bffam=\sevenbf
  \scriptscriptfont\bffam=\fivebf
  \def\bf{\fam\bffam\tenbf}
  \textfont\ttfam=\tentt
  \def\tt{\fam\ttfam\tentt}
  \abovedisplayskip=12pt plus 3pt minus 9pt
  \belowdisplayskip=\abovedisplayskip
  \abovedisplayshortskip=0pt plus 3pt
  \belowdisplayshortskip=4pt plus 3pt 
  \smallskipamount=3pt plus 1pt minus 1pt
  \medskipamount=6pt plus 2pt minus 2pt
  \bigskipamount=12pt plus 4pt minus 4pt
  \normalbaselineskip=12pt
  \setbox\strutbox=\hbox{\vrule height8.5pt depth3.5pt width0pt}
  \let\bigf@nt=\tenrm       \let\smallf@nt=\sevenrm
  \normalbaselines\rm}

\def\eightpoint{
  \textfont0=\eightrm \scriptfont0=\sixrm \scriptscriptfont0=\fiverm
  \def\rm{\fam\z@\eightrm}
  \textfont1=\eighti  \scriptfont1=\sixi  \scriptscriptfont1=\fivei
  \def\oldstyle{\fam\@ne\eighti}\let\old=\oldstyle
  \textfont2=\eightsy \scriptfont2=\sixsy \scriptscriptfont2=\fivesy
  \textfont\gothfam=\eightgoth \scriptfont\gothfam=\sixgoth
  \scriptscriptfont\gothfam=\fivegoth
  \def\goth{\fam\gothfam\eightgoth}
  
  \textfont\itfam=\eightit
  \def\it{\fam\itfam\eightit}
  \textfont\slfam=\eightsl
  \def\sl{\fam\slfam\eightsl}
  \textfont\bffam=\eightbf \scriptfont\bffam=\sixbf
  \scriptscriptfont\bffam=\fivebf
  \def\bf{\fam\bffam\eightbf}
  \textfont\ttfam=\eighttt
  \def\tt{\fam\ttfam\eighttt}
  \abovedisplayskip=9pt plus 3pt minus 9pt
  \belowdisplayskip=\abovedisplayskip
  \abovedisplayshortskip=0pt plus 3pt
  \belowdisplayshortskip=3pt plus 3pt 
  \smallskipamount=2pt plus 1pt minus 1pt
  \medskipamount=4pt plus 2pt minus 1pt
  \bigskipamount=9pt plus 3pt minus 3pt
  \normalbaselineskip=9pt
  \setbox\strutbox=\hbox{\vrule height7pt depth2pt width0pt}
  \let\bigf@nt=\eightrm     \let\smallf@nt=\sixrm
  \normalbaselines\rm}

\tenpoint

\def\pc#1{\bigf@nt#1\smallf@nt}         \def\pd#1 {{\pc#1} }

\catcode`\;=\active
\def;{\relax\ifhmode\ifdim\lastskip>\z@\unskip\fi
\kern\fontdimen2  -1.2 \fontdimen3 \string;}

\catcode`\:=\active
\def:{\relax\ifhmode\ifdim\lastskip>\z@\unskip\fi\penalty\@M\ \fi\string:}

\catcode`\!=\active
\def!{\relax\ifhmode\ifdim\lastskip>\z@
\unskip\fi\kern\fontdimen2  -1.1 \fontdimen3 \string!}

\catcode`\?=\active
\def?{\relax\ifhmode\ifdim\lastskip>\z@
\unskip\fi\kern\fontdimen2  -1.1 \fontdimen3 \string?}

\frenchspacing

\def\raggedbottom{\topskip 10pt plus 36pt\r@ggedbottomtrue}

\def\pointir{\unskip . --- \ignorespaces}

\def\Medbreak{\vskip-\lastskip\medbreak}

\long\def\th#1 #2\enonce#3\endth{
   \Medbreak\noindent
   {\pc#1} {#2\unskip}\pointir{\it #3}\smallskip}

\def\decale#1{\smallbreak\hskip 28pt\llap{#1}\kern 5pt}
\def\decaledecale#1{\smallbreak\hskip 34pt\llap{#1}\kern 5pt}
\def\puce{\smallbreak\hskip 6pt{$\scriptstyle\bullet$}\kern 5pt}

\def\eqalign#1{\null\,\vcenter{\openup\jot\m@th\ialign{
\strut\hfil$\displaystyle{##}$&$\displaystyle{{}##}$\hfil
&&\quad\strut\hfil$\displaystyle{##}$&$\displaystyle{{}##}$\hfil
\crcr#1\crcr}}\,}

\catcode`\@=12

\showboxbreadth=-1  \showboxdepth=-1

\newcount\numerodesection \numerodesection=1
\def\section#1{\bigbreak
 {\bf\number\numerodesection.\ \ #1}\nobreak\medskip
 \advance\numerodesection by1}

\mathcode`A="7041 \mathcode`B="7042 \mathcode`C="7043 \mathcode`D="7044
\mathcode`E="7045 \mathcode`F="7046 \mathcode`G="7047 \mathcode`H="7048
\mathcode`I="7049 \mathcode`J="704A \mathcode`K="704B \mathcode`L="704C
\mathcode`M="704D \mathcode`N="704E \mathcode`O="704F \mathcode`P="7050
\mathcode`Q="7051 \mathcode`R="7052 \mathcode`S="7053 \mathcode`T="7054
\mathcode`U="7055 \mathcode`V="7056 \mathcode`W="7057 \mathcode`X="7058
\mathcode`Y="7059 \mathcode`Z="705A

% handling accented characters in plain TeX :

\def\diagram#1{\def\normalbaselines{\baselineskip=0pt\lineskip=5pt}
\matrix{#1}}

\def\vfl#1#2#3{\llap{$\textstyle #1$}
\left\downarrow\vbox to#3{}\right.\rlap{$\textstyle #2$}}

\def\ufl#1#2#3{\llap{$\textstyle #1$}
\left\uparrow\vbox to#3{}\right.\rlap{$\textstyle #2$}}

\def\hfl#1#2#3{\smash{\mathop{\hbox to#3{\rightarrowfill}}\limits
^{\textstyle#1}_{\textstyle#2}}}

\def\ogoth{{\goth o}}

\def\pgoth{{\goth p}}

\def\P{{\bf P}}
\def\Q{{\bf Q}}
\def\Qp{\Q_p}

\def\R{{\bf R}}

\def\N{{\bf N}}

\def\Z{{\bf Z}}

\def\F{{\bf F}}
\def\Fp{{\F_{\!p}}}

\def\Aut{\mathop{\rm Aut}\nolimits}
\def\Hom{\mathop{\rm Hom}\nolimits}

\def\Int{\mathop{\rm Int}\nolimits}

\def\Card{\mathop{\rm Card}\nolimits}
\def\car{\mathop{\rm car}\nolimits}
\def\Gal{\mathop{\rm Gal}\nolimits}
\def\Ker{\mathop{\rm Ker}\nolimits}

\def\Im{\mathop{\rm Im}\nolimits}
\def\Res{\mathop{\rm Res}\nolimits}
\def\Cor{\mathop{\rm Cor}\nolimits}

\def\droite#1{\,\hfl{#1}{}{8mm}\,}
\def\series#1{(\!(#1)\!)}

\def\to{\rightarrow}

\def\normressym(#1,#2)_#3{\displaystyle\left({#1,#2\over#3}\right)}

\def\mod{\mathop{\rm mod.}\nolimits}
\def\pmod#1{\;(\mod#1)}

\newcount\refno 
\long\def\ref#1:#2<#3>{                                        
\global\advance\refno by1\par\noindent                              
\llap{[{\bf\number\refno}]\ }{#1} \pointir{\it #2} #3\goodbreak }

\def\citer#1(#2){[{\bf\number#1}\if#2\empty\relax\else,\ {#2}\fi]}

\newbox\bibbox
\setbox\bibbox\vbox{\bigbreak
\centerline{{\pc BIBLIOGRAPHY}}

\ref{\pc ARTIN} (E):
Galois Theory,
<University of Notre Dame, Notre Dame, 1942, i+70 pp.>
\newcount\artin \global\artin=\refno

\ref{\pc BROWN} (K):
Cohomology of groups,
<Springer-Verlag, New York, 1994, x+306 pp.>
\newcount\brown \global\brown=\refno

\ref{\pc CASSELS} (J):
Local fields,
<Cambridge University Press, Cambridge, 1986, xiv+360 pp.>
\newcount\cassels \global\cassels=\refno

\ref{\pc DALAWAT} (C):
Local discriminants, kummerian extensions, and elliptic curves,
<J.\ Ramanujan Math.\ Soc. {\bf 25} (2010) 1,
pp.~25--80. Cf.~0711.3878v2.>      
\newcount\locdisc \global\locdisc=\refno

\ref{\pc DALAWAT} (C):
A first course in Local arithmetic,
<arXiv\string:0903.2615v1.>    
\newcount\course \global\course=\refno

\ref{\pc DALAWAT} (C):
Further remarks on local discriminants, 
<J.\ Ramanujan Math.\ Soc.\ {\bf 25} (2010) 4, pp.~391--417. 
Cf.~arXiv\string:0909.2541v1.>  
\newcount\further \global\further=\refno

\ref{\pc DALAWAT} (C):
Final remarks on local discriminants, 
<J.\ Ramanujan Math.\ Soc.\ {\bf 25} (2010 ) 4, pp.~419--432. 
Cf.~arXiv\string:0912.2829v2.> 
\newcount\final \global\final=\refno

\ref{\pc DEL \pc CORSO} (I) and {\pc DVORNICICH} (R):
The compositum of wild extensions of local fields of prime degree,
<Monatsh.\ Math.\ {\bf 150} (2007) 4, pp.~271--288.>
\newcount\delcorso \global\delcorso=\refno

\ref{\pc DOUD} (D):
Wild ramification in number field extensions of prime degree,
<Arch.\ Math.\ (Basel) {\bf 81} (2003) 6, pp.~646--649.>
\newcount\doud \global\doud=\refno

\ref{\pc HASSE} (H):
Zahlentheorie,
<Akademie-Verlag, Berlin, 1969, 611~pp.>
\newcount\hasse \global\hasse=\refno

\ref{\pc HUPPERT} (B):
Endliche Gruppen. I. 
<Springer-Verlag, Berlin-New York 1967, xii+793 pp.>
\newcount\huppert \global\huppert=\refno

% \ref{\pc KRASNER} (M): 
% Remarques au sujet d'une note de J.-P.~Serre: ``Une `formule de masse' pour
% les extensions totalement ramifi{\'e}es de degr{\'e} donn{\'e} d'un corps
% local'' : une d{\'e}monstration de la formule de M.~Serre {\`a} partir de mon
% th{\'e}or{\`e}me sur le nombre des extensions s{\'e}parables d'un corps
% valu{\'e} localement compact, qui sont d'un degr{\'e} et d'une diff{\'e}rente
% donn{\'e}s, 
% <Comptes Rendus {\bf 288} (1979)~18, pp.~A863--A865.>
% \newcount\krasner \global\krasner=\refno

\ref{\pc LENSTRA} (H):
A normal basis theorem for infinite Galois extensions,
<Nederl.\ Akad.\ Wetensch.\ Indag.\ Math.\ {\bf 47} (1985)~2,
pp.~221--228.>
\newcount\lenstra \global\lenstra=\refno

\ref{\pc NEUMANN} (O):
Two proofs of the Kronecker-Weber theorem ``according to Kronecker, and
Weber'',
<J.\ Reine Angew.\ Math.\  {\bf 323}  (1981), pp.~105--126.>
\newcount\neumann \global\neumann=\refno

% \ref{\pc ROBINSON} (D):
% A course in the theory of groups,
% <Springer-Verlag, New York-Berlin, 1982, xvii+481~pp.>
% \newcount\robinson \global\robinson=\refno

% \ref{\pc ROTMAN} (J):
% An introduction to the theory of groups,
% <Springer-Verlag, New York, 1995, xvi+513~pp.>
% \newcount\rotman \global\rotman=\refno

\ref{\pc SERRE} (J-P):
Corps locaux,
<Publications de l'Universit{\'e} de Nancago {\sevenrm VIII}, Hermann,
Paris, 1968, 245 pp.>
\newcount\corpslocaux \global\corpslocaux=\refno

\ref{\pc SERRE} (J-P):
Une ``\thinspace formule de masse$\,$" pour les extensions totalement
ramifi{\'e}es de degr{\'e} donn{\'e} d'un corps local, 
<Comptes Rendus {\bf 286} (1978), pp.~1031--1036.>
\newcount\serremass \global\serremass=\refno

} %\bibbox

\centerline{\bf Serre's {\it ``\thinspace formule de masse$\,$''\/} in
  prime   degree}  
\bigskip\bigskip 
\centerline{Chandan Singh Dalawat} 
\centerline{Harish-Chandra Research Institute}
\centerline{Chhatnag Road, Jhunsi, Allahabad 211019, India} 
\centerline{dalawat@gmail.com}

\bigskip\bigskip

% {\it Dans la th{\'e}orie des {\'e}quations, j'ai cherch{\'e} dans quel cas les
%   {\'e}quations {\'e}tait r{\'e}solubles par des radicaux.}\hfill ---
%   {\'E}variste   Galois, 29 May 1932.
% 
% \bigskip\bigskip

{{\bf Abstract}.  For a local field $F$ with finite residue field of
  characteristic~$p$, we describe completely the structure of the
  filtered $\Fp[G]$-module $K^\times\!/K^{\times p}$ in
  characteristic~$0$ and $K^+\!/\wp(K^+)$ in characteristic~$p$, where
  $K=F(\!\root{p-1}\of{F^\times})$ and $G=\Gal(K|F)$.  As an
  application, we give an elementary proof of Serre's mass formula in
  degree~$p$.  We also determine the compositum $C$ of all degree-$p$
  separable extensions with solvable galoisian closure over an
  arbitrary base field, and show that $C$ is $K(\!\root
  p\of{K^\times})$ or $K(\wp^{-1}(K))$ respectively, in the case of
  the local field $F$.  Our method allows us to compute the
  contribution of each character $G\to\F_p^\times$ to the degree-$p$
  mass formula, and, for any given group $\Gamma$, the contribution of
  those degree-$p$ separable extensions of $F$ whose galoisian closure
  has group $\Gamma$.  \footnote{}{Keywords~: {\it Formule de masse de Serre,
  Serre'schen Ma\ss formel}, Serre's mass formula.}}

\bigskip

{\bf 1. Introduction}\pointir Let $p$ be a prime number and let $F$ be a local
field with finite residue field $k$ of characteristic~$p$ and cardinality
$q=p^f$, so that the characteristic $\car(F)$ of $F$ either~$p$ or~$0$ and
indeed $F=k(\!(\pi)\!)$ in the characteristic-$p$ case and $F$ is a finite
extension of $\Qp$ of residual degree~$f$ and ramification index $e=[F:\Qp]/f$
in the characteristic-$0$ case~; put $e=+\infty$ in the former case.

For a totally ramified separable extension $E$ of $F$ of degree~$n$ and
discriminant $\delta_{E|F}$ of valuation $v(\delta_{E|F})$, put
$$
c(E)=v(\delta_{E|F})-(n-1),
$$
so that $c(E)=0$ if and only if $E|F$ is tamely ramified or equivalently
$n$~is prime to~$p$.  Serre's mass formula says that when $E$ runs through
such extensions (contained in a fixed separable algebraic closure of $F$),
then
$$
\sum_E q^{-c(E)}
=n.
$$
This beautiful formula is easy to prove in the tame case when $n$ is prime
to~$p$ but lies much deeper when $p|n$.  The main purpose of this Note is to
give an elementary proof of this formula in the case $n=p$ and to develop the
algebraic and arithmetic ingredients on which this proof relies.

(It turns out that the method adopted here leads to several
refinements of the mass formula.  For instance, one can compute the
contribution of those $E$ which are cyclic over $F$, or those $E$
whose galoisian closure $\tilde E$ over $F$ is of the form $EF'$ for
some unramified extension $F'|F$ (depending on $E$), or those $E$ for
which the group $\Gal(\tilde E|F)$ is isomorphic to a given group
$\Gamma$.  See \S9 for details.)

There are two sources of inspiration for our method.  The first one is our
recent reworking \citer\final() of a standard technique used in proofs of the
local \citer\cassels(p.~155) or the global \citer\neumann(p.~110)
Kronecker-Weber theorem.  It allowed us to compute the contribution of
degree-$p$ {\it cyclic\/} extensions to the mass formula \citer\final(),~\S6.

The idea was that in characteristic~$p$, the set of such extensions is in
bijection with the set of $\Fp$-lines in $F/\wp(F)$, where $\wp:x\mapsto
x^p-x$, and that in characteristic~$0$ it is in bijection with the set of
$\Fp$-lines in the $\omega$-eigenspace for the action of $\Delta$ on
$F(\zeta)^\times\!/F(\zeta)^{\times p}$, where $\zeta$ is a primitive $p$-th
root of~$1$, $\Delta=\Gal(F(\zeta)|F)$, and $\omega:\Delta\to\F_p^\times$ is
the cyclotomic character giving the action of $\Delta$ on the $p$-th roots
of~$1$.  Roughly speaking, what we did for $\omega$ there, we need to do here
for all characters $\chi:G\to\F_p^\times$, where $G=\Gal(K|F)$ and
$K=F(\!\root{p-1}\of{F^\times})$.

The results and methods of Del Corso and Dvornicich \citer\delcorso() provided
the second source of inspiration and ideas. They study the action of $G$ on
$K^\times\!/K^{\times p}$ in characteristic~$0$ and prove that the compositum
of all \hbox{degree-$p$} extensions of $F$ is $K(\!\root p\of{K^\times})$. The
\hbox{characteristic-$p$} analogue of their main theorem, stating that the
compositum of all \hbox{degree-$p$} separable extensions of $F$ is
$K(\wp^{-1}(K))$, where $K$ is still $F(\!\root{p-1}\of{F^\times})$, can be
found below (prop.~36).  We also study (\S8) the compositum of all degree-$p'$
extensions of $F$ for a prime $p'\neq p$, which turns out to be $K'(\!\root
p'\of{K'^\times})$, where $K'=F(\zeta')$ and $\zeta'$ is a primitive $p'$-th
root of~$1$.

Our first task is thus to extend the methods and the main result of
\citer\delcorso() to characteristic-$p$ local fields.  They turn out to behave
exactly as their \hbox{characteristic-$0$} counterparts would have if
$e=+\infty$ and $\omega=1$.  We give a unified presentation of the two cases
which is perhaps more intrinsic and conceptual than the treatment of the
characteristic-$0$ case in \citer\delcorso(), to which our debt should however
be clear to the reader.

Our main arithmetic contribution is thus an explicit description of the
structure of the filtered $G$-module $K^\times\!/K^{\times p}$
(resp.~$K^+\!/\wp(K^+)$).  It is $\Fp[G]$-isomorphic to
$$
\Fp\{\omega\}\oplus k[G]^e\oplus\Fp\qquad 
(\hbox{resp. }\Fp\oplus k[G]\oplus k[G]\oplus\cdots)
$$
(where $\Fp\{\omega\}$ denotes an $\Fp$-line on which $G$ acts via
$\omega$) with a specific filtration to be described in detail in \S6.  For
degree-$p'$ extensions, the analogue is the filtered $\F_{p'}[G']$-module
$\F_{p'}\{\omega'\}\oplus\F_{p'}$, where $\omega'$ is the mod-$p'$ cyclotomic
character of $G'=\Gal(K'|F)$~: it is as if $e=0$.  We then use this
description to give an elementary proof of Serre's mass formula in prime
degree~; it is uniformly applicable to all three cases~: $0<e<+\infty$,
$e=+\infty$, and $e=0$.

% equicharacteristic ($0<e<+\infty$), mixed-characteristic ($e=+\infty$), and
% $p'\neq p$ ($e=0$).
% 
Let us give a sketch of the proof in the degree-$p$ case. Every separable
degree-$p$ extension $E$ of $F$ becomes cyclic when translated to $K$.  It
corresponds therefore to an $\Fp$-line $D$ in $K^\times\!/K^{\times p}$
(resp.~$K^+\!/\wp(K^+)$).  Such lines are stable under the action of~$G$, and every
$G$-stable line arises from some $E$.  Two such extensions $E,E'$ give rise to
the same~$D$ if and only if they are conjugate over~$F$.  If $E|F$ is not
cyclic, then it has $p$ conjugates.

The invariant $c(E)$ of $E$ can be recovered from the ``\thinspace
level\thinspace'' of~$D$ in the filtration on $K^\times\!/K^{\times p}$
(resp.~$K^+\!/\wp(K^+)$) --- the integer $m$ such that $D\subset\bar U_m$ but
$D\not\subset\bar U_{m+1}$ (resp.~$D\subset\overline{\pgoth^m}$ but
$D\not\subset\overline{\pgoth^{m+1}}$), where $(\bar U_n)_{n\in\N}$
(resp.~$(\overline{\pgoth^n})_{n\in\Z}$)is the induced filtration --- by the
{\it Schachtelungssatz}.  The number $r$ of $E$ giving rise to $D$ can be read
off from the character $\chi:G\to\F_p^\times$ through which $G$ acts on $D$~;
we have $r=1$ if $\chi=\omega$ and $r=p$ if $\chi\neq\omega$, where $\omega$
is the mod-$p$ cyclotomic character in characteristic~$0$ and $\omega=1$ by
convention in characteristic~$p$.  So the sum over all ramified separable
degree-$p$ extensions $E$ of $F$ gets replaced by a sum over all $G$-stable
lines $D$ in $K^\times\!/K^{\times p}$ or $K^+\!/\wp(K^+)$ other than the line $\bar
U_{pe}$ in characteristic~$0$ and other than the line $\bar\ogoth$ in
characteristic~$p$.  As we have already determined (\S6) the structure of the
filtered $\Fp[G]$-module $K^\times\!/K^{\times p}$ (resp.~$K^+\!/\wp(K^+)$), the
contribution from these $D$ can be computed level by level, leading to the
result.  See \S7 for the complete proof and \S9 for an extended
summary along with some refinements.

This proof may be thought of a generalisation from the case $p=2$
(\citer\locdisc(prop.~67), \citer\final(prop.~14)) where it amounted to a
trivial identity because every separable quadratic extension is kummerian ---
of the form $F(x)$ with $x^2\in F^\times$ (resp.~$x^2-x\in F$).  The reader
may also wish to contrast this proof with the much easier case of degree-$p'$
extensions for a prime $p'\neq p$ outlined in \S8.

\medskip
\line{\hfill\hbox{***}\hfill}
\medbreak

While this application to the degree-$p$ case of Serre's mass formula (\S7)
for local fields of residual characteristic~$p$ was our main motivation, the
contribution of this Note is not confined to this proof.  We have first
developed (\S4) the purely algebraic theory of solvable extensions of prime
degree (def.~11).  Observations such as {\it Remark\/}~13 and lemmas~17--20
are of independent interest.  The arithmetic theory of filtered galoisian
modules (\S\S6,8) has been developed to the fullest (in the case at hand) in an
intrinsic manner without distinction between the mixed-characteristic, the
equicharacteristic, and the $p'\neq p$ cases.  

The general algebraic theory (\S4) and the local arithmetic theory
(\S\S6,8) have other applications.  For instance, the degree-$p$ mass
formula can be thought of as the sum of the contributions of various
characters $\xi:G\to\F_p^\times$.  Our method lets us compute
the contribution of each $\xi$, which allows us to compute
the contribution of various natural classes of degree-$p$ separable
extensios of $F$  (see \S9).

\bigbreak

{\bf 2. Contents}\pointir \S3 recalls some basic facts about groups and their
representations which are applied in \S4 to study degree-$p$ solvable
extensions (def.~11) of an arbitrary base field, culminating in the
determination of the compositum of all such extensions (lemma~20).  What is
needed for the subsequent \S\S\ is summarised at the end of \S4.

We return to our local field $F$ in \S5 and work out the consequences
of the algebraic theory in this particular case.  We undertake in \S6
a detailed study of the filtered $G$-module $K^\times\!/K^{\times p}$
(resp.~$K^+\!/\wp(K^+)$) of $\Fp$-dimension $2+(p-1)^2ef$ (resp.~$+\infty$),
where $G=\Gal(K|F)$ and $K=F(\!\root p-1\of{F^\times})$, leading to
the determination of the compositum of all separable extensions of $F$
of degree~$p$ (prop.~36), and to the proof of Serre's degree-$p$ mass
formula (th.~39) in \S7.  The method is supple enough to compute the
contributions of individual characters $\chi\in F^\times\!/F^{\times
p-1}$, as explained in \S9.  Finally, we indicate (\S8) how simple the
whole theory becomes when dealing with degree-$p'$ extensions for some
prime $p'\neq p$ ($K$ gets replaced by $K'=F(\!\root p'\of1)$ and the
$\F_{p'}$-space $K'^\times\!/K'^{\times p'}$ is $2$-dimensional).

\bigbreak

{\bf 3. Groups and their representations}\pointir All we need are some fairly
standard results but we include some proofs.

First, let $P\subset{\goth S}_p$ be a subgroup of order~$p$ and $N$ its
normaliser~; identify $\Aut P$ with $\F_p^\times$.  For every $\eta\in N$,
conjugation by $\eta$ induces an automorphism $\Int(\eta)\in\F_p^\times$ of
$P$ which depends only on the class $\bar\eta\in N/P$, leading to a map
$N/P\to\F_p^\times$.

\th LEMMA 1
\enonce
The map\/ $N/P\to\F_p^\times$ is an isomorphism~: we have a split short exact
sequence\/ $1\to P\to N\to\F_p^\times\to1$.
\endth
In other words, $N=P\times_\iota\F_p^\times$, where
$\iota:a\mapsto(\sigma\mapsto\sigma^a)$ is the canonical isomorphism 
$\F_p^\times\to\Aut P$.  Consequently, a subgroup $\Gamma\subset{\goth
  S}_p$ containing $P$ and contained in $N$ is transitive and solvable, and
indeed $\Gamma$ is the canonical split extension of a subgroup of
$\F_p^\times$ by $P$. There is a converse.   

\th LEMMA 2 (Galois)
\enonce
A transitive subgroup\/ $\Gamma\subset{\goth S}_p$ is solvable if and only if
it contains a unique Sylow $p$-subgroup\/ $P$ (of order\/~$p$).
\endth
This was guessed at from the reference to \citer\artin(ch.~3, th.~7) in
\citer\delcorso()~; Robin Chapman and Jack Schmidt provided short arguments
on MathOverflow, reproduced below, and Matthew Emerton indicated the
provenance. 

The order of a transitive subgroup $\Gamma\subset{\goth S}_p$ is divisible
by~$p$ (because the orbit has $p$ elements) but not by $p^2$ (because the
order $p!$ of ${\goth S}_p$ is not). Therefore $\Gamma$ has a Sylow
$p$-subgroup $P$ of order~$p$.  If it is unique, then $P\subset\Gamma$ is
normal, and hence $\Gamma$ is contained in the normaliser $N$ of $P$ in
${\goth S}_p$. By lemma~1, $\Gamma$ is solvable and the result follows in this
case.  We show next that this is the only case~: $P$ is unique if $\Gamma$ is
solvable.
  
Suppose that $\Gamma$ has Sylow $p$-subgroups other than~$P$~; they have to be
conjugates of $P$ in $\Gamma$ (Sylow). If $H\subset\Gamma$ is a normal
subgroup of order $>1$, then $H$ must be transitive since otherwise the orbits
of $H$ would form a partition of $\{1,\cdots,p\}$ invariant under the action
of $\Gamma$, which is impossible because $p$~is prime. So $H$ contains some
(and hence every) conjugate of $P$.  It follows that $\Gamma$ is not solvable,
as there are no normal subgroups containing only one conjugate of~$P$.  A
standard reference is \citer\huppert(p.~163).

\smallskip

Next, let $1\to C\to\Gamma\to G\to1$ be any short exact sequence of finite
groups in which $C$ is {\it commutative}.

\th LEMMA 3 (Schur-Zassenhaus)
\enonce
If the orders of\/ $C$ and $G$ are mutually prime, then 
sections\/ $G\to\Gamma$ do exist, and any two sections\/ $G\to\Gamma$differ
by an inner automorphism\/ $\Int(\gamma):\Gamma\to\Gamma$ for some\/
$\gamma\in C$.  
\endth
Associating to the image $\bar\gamma\in G$ of $\gamma\in\Gamma$ the
automorphism $\sigma\mapsto\gamma\sigma\gamma^{-1}$ of $C$ (which depends only
on $\bar\gamma$ because $C$ is commutative), we get an action $\theta:G\to\Aut
C$ of $G$ on $C$.  The group $H^2(G,C)_\theta$ (the action $\theta$ is often
omitted from the notation) classifies extensions of $G$ by $C$ in which $G$
acts on $C$ via~$\theta$~; the twisted product $C\times_\theta G$ is the
neutral element of $H^2(G,C)_\theta$.  For every subgroup $G'\subset G$, we
have the natural maps
$$
\Res:H^2(G,C)_\theta\to H^2(G',C)_\theta,\quad
\Cor:H^2(G',C)_\theta\to H^2(G,C)_\theta,
$$
of restriction and corestriction whose composite $\Cor\circ\Res$ is
multiplication by the index $\hbox{\bf(}G:G'\hbox{\bf)}$.  Upon taking
$G'=\{1\}$, we see that $H^2(G,C)_\theta$ is killed by the order of $G$.  But
it is also killed by the order of ~$C$. As these orders are mutually prime by
hypothesis, we have $H^2(G,C)_\theta=0$, and hence $\Gamma=C\times_\theta G$
is a twisted product, which admits sections.

Now, $C$-conjugacy classes of sections $G\to C\times_\theta G$ are classified
by the group $H^1(G,C)_\theta$, which can be shown to vanish exactly as above.
I thank Joseph Oesterl{\'e} for supplying this argument at short notice.
See \citer\brown() for details and \citer\corpslocaux(IV, 2, cor.~4) for a
direct proof in the case needed below, in which $C$ is a $p$-group and $G$ is
commutative of exponent dividing $p-1$.

{\it Remark 4}\pointir When moreover $G$ is commutative and the action $\theta$
is trivial, then there is a unique section $G\to C\times G$, namely the
inclusion. 

{\it Remark 5}\pointir A solvable transitive subgroup $\Gamma\subset{\goth
  S}_p$ has exactly $p$ subgroups of index~$p$, unless $\Gamma$ is commutative
(in which case it is in fact cyclic of order~$p$).  Indeed, $P$ being the
unique Sylow $p$-subgroup of $\Gamma$ (lemma~2), $\Gamma/P$ is a subgroup of
$\F_p^\times$ (lemma~1), hence there are $p$ sections $\Gamma/P\to\Gamma$
unless $\Gamma=P$ (lemma~3).

{\it Remark 6}\pointir Let $\bar\tau\in\F_p^\times$ be a generator of
$\Gamma/P$ and $g$ its order (so that $g\,|\,p-1$).  The groups $P$ and
$\Gamma$ admit the presentations
$$
P=\langle\sigma\;|\;\sigma^p=1\rangle,\quad
\Gamma=\langle\sigma,\tau\;|\;
\sigma^p=1,\tau^g=1,\tau\sigma\tau^{-1}=\sigma^{\bar\tau}\rangle, 
$$
and we have $\Gamma=P\Leftrightarrow\bar\tau=1\Leftrightarrow g=1$.  If
$\Gamma\neq P$, then the $p$ index-$p$ subgroups $G_i\subset\Gamma$ are
generated respectively by $\sigma^i\tau\sigma^{-i}$ as $i$ runs through $\Fp$.
For $i\neq j$, we have $G_i\cap G_j=\{1\}$ and $G_i\cup G_j$ generates
$\Gamma$.

{\it Remark 7}\pointir A version of lemma~3 remains valid for profinite
groups.  We omit the details \footnote{(*)}{(2011/03/15) See for example
Iwasawa ({\it Transactions AMS}, {\bf 80} (1955), 448--469, lemma~5)}, for
$\Gamma$ is finite in our main application to local fields.  The profinite
version would be convenient for the proof of lemma~16 where $1\to
C\to\Gamma\to G\to1$ comes from a galoisian tower $L|K|F$ in which $L|K$ is
cyclic of degree~$p$ and $K|F$ is abelian of exponent dividing $p-1$ (but of
possibly infinite degree).  Writing $L=K(x)$, taking $K'$ to be a finite
galoisian extension of $F$ in $K$ containing the coefficients of the minimal
polynomial $f\in K[T]$ of $x$ and such that $K'(x)$ is cyclic (of degree~$p$)
over $K'$ and galoisian over $F$, lemma~3 as stated above is applicable to the
short exact sequence associated to the galoisian tower $K'(x)|K'|F$~; this
will suffice.

\smallskip

Finally, we need a few facts about representations of a commutative group $G$
of exponent $d$ over a field $k$ of characteristic prime to the order of $G$
and containing all $d$-th roots of~$1$.  So let $V$ be a $k$-space on which
$G$ acts by $k$-automorphisms~; in other words, $V$ is a $k[G]$-module.  For
every character $\chi:G\to k^\times$, we denote by $V(\chi)$ the
$\chi$-eigenspace for the action of $G$ on $V$~; it consists of all $x\in V$
such that $\sigma.x=\chi(\sigma)x$ for every $\sigma\in G$.

\th LEMMA 8
\enonce
The $k$-space\/ $V$ is the internal direct sum of its subspaces\/ $V(\chi)$,
indexed by all characters\/ $\chi:G\to k^\times$.  Every\/ $G$-stable subspace
$W\subset V$ has a\/ $G$-stable supplement\/ $W'$, and the canonical map\/
$W'\to V/W$ is an isomorphism of\/ $k[G]$-modules.   
\endth
Indeed, $k[G]$ is a semisimple algebra and the only simple $k[G]$-modules are
$k$-lines on which $G$ acts through some character $G\to k^\times$.

A particularly interesting example occurs when $G=\Gal(l|k)$ for $l|k$ a
finite abelian extension of exponent $d$ prime to the characteristic of $k$
and such that $k$ contains a primitive $d$-th root of~$1$, and $V=l$.  The
normal basis theorem --- the $k[G]$-module $l$ is free of rank~$1$ --- implies

\th LEMMA 9
\enonce
For every character\/ $\chi:G\to k^\times$ of\/ $G$, the $k$-space\/ $l(\chi)$
is of dimension\/~$1$.
\endth
A direct proof in the special case we need can be found in prop.~30.

{\it Remark 10}\pointir For a version of the normal basis theorem applicable to
$l|k$ of possibly infinite degree, see \citer\lenstra().  The ring $k[G]$ is
replaced by $k[[G]]$ --- the inverse limit of $k[G/H]$ as $H$ runs through open
normal subgroups of $G$, and $l$ is replaced by the inverse limit of $l^H$
under the trace maps~; as a $k[[G]]$-module, the latter limit is free of
rank~$1$.  %We won't need this generalisation.

\bigbreak 

{\bf 4.  Fields and their extensions}\pointir We work over an arbitrary
(commutative) field $F$ whose characteristic $\car(F)$ may or may not be equal
to our fixed prime ~$p$.  If $\car(F)\neq p$, denote by $\omega$ the mod-$p$
cyclotomic character (giving the action on the $p$-th roots of~$1$) and take
$\omega=1$ to be the trivial character if $\car(F)=p$, so that $\omega$ has
values in $\F_p^\times$.

\th DEFINITION 11
\enonce
A degree-$p$ extension\/ $E$ of\/ $F$ is called solvable if it is separable
and if the group\/ $\Gal(\tilde E|F)$ of\/ $F$-automorphisms of its galoisian
closure\/ $\tilde E$ over\/ $F$ is solvable. 
\endth
The terminology is not standard (unless $\tilde E=E$) but unlikely to
confuse. 

\th LEMMA 12
\enonce
For every degree-$p$ solvable extension\/ $E|F$, there exists a
cyclic extension\/ $F'|F$ of degree dividing\/~$p-1$ such that\/ $EF'$ is
cyclic (of degree\/~$p$) over $F'$ and galoisian over\/ $F$.  If\/ $E|F$ is
not cyclic, then it has exactly\/ $p$ conjugates over\/ $F$.
\endth 
Let $\tilde E$ be the galoisian closure of $E$ and $\Gamma=\Gal(\tilde E|F)$~;
the group $\Gamma$ is solvable (def.~11).  Also, $\Gamma$ operates transitively
on the set of $F$-embeddings of $E$ (in any fixed but arbitrary separable
algebraic closure of $F$), so $\Gamma\subset{\goth S}_p$.

Lemma~2 then furnishes a (unique) order-$p$ subgroup $P\subset \Gamma$.  It is
then clear that we may take $F'={\tilde E}^P$.  Indeed, $F'|F$ is cyclic of
group $\Gamma/P$ of order dividing\/~$p-1$ (lemma~1), and $EF'=\tilde E$ is
cyclic over $F'$ and galoisian over $F$.  Finally, {\it Remark\/}~5 says that
if the group $\Gamma$ is not cyclic, then it has $p$ index-$p$ subgroups, so
the extension $E$ has precisely $p$ conjugates over $F$ if it is not cyclic
--- something which is also otherwise clear.  

{\it Remark~13}\pointir If $E'$, $E''$ are two distinct conjugates of $E$,
then $E'E''=\tilde E$ ({\it Remark\/}~6).  Consequently, if a solvable
separable irreducible degree-$p$ polynomial over $F$ has two roots in an
extension $R$ of $F$, then all its roots are in $R$.  The case $R=\R$
\citer\artin(p.~67) is sometimes attributed to Kronecker~: ({\it Wenn eine
  irreductible Gleichung mit ganzzahligen Co{\"e}fficienten aufl{\"o}sbar und
  der Grad derselben eine ungrade\/ ({\rm sic}) Primzahl ist, so sind
  entweder\/ {\rm alle} ihre Wurzeln oder nur\/ {\rm eine} reell} (Werke, IV,
p.~25), but, as Kronecker himself points out, Galois had the general result~:
{\it Pour qu'une {\'e}quation de degr{\'e} premier soit r{\'e}soluble par
  radicaux, il faut et il suffit que deux quelconques de ces racines {\'e}tant
  connues, les autres s'en d{\'e}duisent rationnellement} (Bulletin de M.\ 
F{\'e}russac, XIII (avril 1830), p.~271).  Note that Kronecker's observation
applies to every $R$ because Galois's observation applies to every solvable
transitive group of degree~$p$.

Lemma 12 admits a converse~:

\th LEMMA 14
\enonce
Let\/ $F'|F$ be cyclic of degree dividing\/ $p-1$, and\/ $L|F'$ cyclic of
degree\/~$p$.  If\/ $L|F$ is galoisian, then there exists a (solvable)
degree-$p$ extension\/ $E|F$ such that\/ $L=EF'$, any two such extensions are
conjugate over\/ $F$, and every conjugate of\/ $E$ is contained in\/ $L$.
\endth
Suppose that $L|F$ is galoisian of group $\Gamma=\Gal(L|F)$.  We then have an
exact sequence $1\to P\to\Gamma\to G\to1$, with $P=\Gal(L|F')$ of order~$p$
and $G=\Gal(F'|F)$ (cyclic) of order dividing~$p-1$.  The Schur-Zassenhaus
theorem (lemma~3) then implies that $\Gamma$ has a subgroup $G'$ of index~$p$,
and that any two such subgroups are conjugate in $\Gamma$.  It follows that
$E=L^{G'}$ is a solvable (def.~11) degree-$p$ extension of $F$ in $L$ such
that $L=EF'$, and that any two such extensions are conjugate over~$F$.
Finally, every $F$-conjugate of $E$ is contained in $L$ because $L|F$ is
galoisian and contains~$E$.

\smallskip

Let $K$ be the compositum of all cyclic extensions of $F$ of degree
dividing~$p-1$, so that $K$ is the maximal abelian extension of $F$ of
exponent dividing $p-1$.

\th LEMMA 15
\enonce
For every degree-$p$ solvable extension\/ $E$ of\/ $F$, the compositum
$EK$ is cyclic over $K$ and galoisian over\/ $F$. 
\endth 
This follows from lemma~12, the definition of $K$, and def.~$11$.

\smallskip

Which cyclic extensions $L|K$ arise as $L=EK$ for some (degree-$p$,
solvable) extension $E|F$~?  If $L$ arises from $E$, then $L$ would be
galoisian over~$F$ (lemma~15).  Conversely,

\th LEMMA 16
\enonce
If a degree-$p$ cyclic extension\/ $L$ of\/ $K$ is galoisian over\/ $F$, then 
there is a degree-$p$ solvable extension\/ $E|F$ such that\/ $L=EK$~; two
such extensions $E,E'$ give rise to the same $L$ if and only if they are
conjugate over\/ $F$, and every conjugate of\/ $E$ is contained in\/ $L$.
\endth
This follows from the Schur-Zassenhaus theorem (lemma~3) exactly in the same
way as lemma~14 does (cf.~{\it Remark\/}~7).

\smallskip

But the great thing about degree-$p$ cyclic extensions $L|K$ is that they
correspond bijectively to lines $D$ in the $\Fp$-space $K^\times\!/K^{\times
  p}$ in case the characteristic of $F$ is $\neq p$ because $K^\times$
contains a primitive $p$-th root of~$1$, and in the space $K^+\!/\wp(K^+)$ in the
characteristic-$p$ case.  When is the (degree-$p$, cyclic) extension of $K$
corresponding to $D$ galoisian over $F$~?  Precisely when $D$ is stable under
the $G$-action on these spaces, where $G=\Gal(K|F)$.

\th LEMMA 17
\enonce
Let\/ $D$ be a line in\/ $K^\times\!/K^{\times p}$ (resp.~$K^+\!/\wp(K^+)$) and let\/
$L=K(\root p\of D)$ (resp.~$L=K(\wp^{-1}(D)$) be the corresponding cyclic
extension of degree\/~$p$.  Then\/ $L|F$ is galoisian if and only if\/ $D$
is\/  $G$-stable.  If\/ $G$ acts on\/ $D$ via the character\/
$\chi:G\to\F_p^\times$, then it acts on\/ $P=\Gal(L|K)$ via the
character\/ $\omega\chi^{-1}$, where $\omega:G\to\F_p^\times$ is the
cyclotomic character if\/ $\car(F)\neq p$ and\/ $\omega=1$ if\/ $\car(F)=p$.
\endth
Note first of all that when $L|F$ is galoisian of group $\Gamma=\Gal(L|F)$, we
have a short exact sequence $1\to P\to\Gamma\to G\to1$ which provides an
action of $G$ on $P$ by conjugation and hence a character
$\xi:G\to\F_p^\times$.  It is being asserted that $\xi=\omega\chi^{-1}$.  Note
also the corollary that $L|F$ is abelian if and only if $G$ acts on $D$ via
$\omega$.

We have $D=\Ker(\iota:K^\times\!/K^{\times p}\to L^\times\!/L^{\times p})$ in
characteristic~$\neq p$ and $D=\Ker(\iota:K^+\!/\wp(K^+)\to L^+/\wp(L^+))$ in
characteristic~$p$, where $\iota$ is induced by the inclusion of $K$ in $L$.

\goodbreak

If $L|F$ is galoisian of group $\Gamma$, then
$\gamma\circ\iota=\iota\circ\gamma$ for every $\gamma\in\Gamma$, which shows
that $D=\gamma(D)$ and hence $D$ is $G$-stable.  Conversely, if $D$ is
$G$-stable, then $g(a)\in D$ for every $a\in D$ and every $g\in G$.  Therefore
$L$ contains a $p$-th root (resp.~$\wp$-th root) of $g(a)$ for every $a\in D$
and for every $g\in G$, making it galoisian over~$F$.

Suppose now that $L|F$ is galoisian and that $G$ acts on $D$ via $\chi$.  Let
$1\to P\to\Gamma\to G\to1$ be the associated short exact sequence of groups.
Let us first show the final assertion in characteristic~$\neq p$.  Write
$L=K(x)$ for some $x\in L^\times$ such that $a=x^p$ is in $K^\times$ and $\bar
a\in D$ is a generator.  Also choose a generator $\sigma\in P$ and denote by
$\zeta\in{}_p\mu$ the $p$-th root of~$1$ such that $\sigma(x)=\zeta x$~; we
have to show that
$\gamma\sigma\gamma^{-1}=\sigma^{\omega\chi^{-1}(\bar\gamma)}$ for every
$\gamma\in\Gamma$ of image $\bar\gamma\in G$, for which it is sufficient to
show that $\gamma\sigma\gamma^{-1}(x)=\zeta^{\omega\chi^{-1}(\bar\gamma)}x$.

All that is given to us is that $\bar\gamma(\bar a)={\bar
  a}^{\chi(\bar\gamma)}$ for every $\gamma\in\Gamma$.  As
$\gamma(x)^p=\gamma(x^p)=\gamma(a)$, we must have $\gamma(x)=b_\gamma
x^{\chi(\bar\gamma)}$ for some $b_\gamma\in K^\times$ with
$\gamma^{-1}(b_{\gamma})b_{\gamma^{-1}}^{\chi(\bar\gamma)}=1$ (modulo
$K^{\times p}$), to ensure that $\gamma^{-1}\gamma(x)=x$.  Now, 
% $$
% \gamma\sigma\gamma^{-1}(x)
% =\gamma\sigma(b_{\gamma^{-1}}x^{\chi^{-1}(\bar\gamma)})
% =\gamma(b_{\gamma^{-1}}
%      \zeta^{\chi^{-1}(\bar\gamma)}x^{\chi^{-1}(\bar\gamma)})
% $$
% which equals $\gamma(b_{\gamma^{-1}})
% \zeta^{\omega\chi^{-1}(\bar\gamma)}b_\gamma^{\chi^{-1}(\bar\gamma)}x
% =\zeta^{\omega\chi^{-1}(\bar\gamma)}x$, 
$$
\eqalign{
\gamma\sigma\gamma^{-1}(x)
&=\gamma\sigma(b_{\gamma^{-1}}x^{\chi^{-1}(\bar\gamma)})\cr
&=\gamma(b_{\gamma^{-1}}
     \zeta^{\chi^{-1}(\bar\gamma)}x^{\chi^{-1}(\bar\gamma)})\cr
&=\gamma(b_{\gamma^{-1}}) 
\zeta^{\omega\chi^{-1}(\bar\gamma)}b_\gamma^{\chi^{-1}(\bar\gamma)}x\cr
&=\zeta^{\omega\chi^{-1}(\bar\gamma)}x,\cr
}
$$
at least modulo $K^{\times p}$, proving the result in this case.  A similar
argument works in characteristic~$p$ upon replacing the multiplicative
notation with the additive notation.

Let's do the transcription.  Write $L=K(x)$ for some $x\in L$ such that
$a=x^p-x$ is in $K$ and $\bar a\in D$ is a generator.  Let $\sigma\in P$ be
the generator such that $\sigma(x)=x+1$~; we
have to show that
$\gamma\sigma\gamma^{-1}=\sigma^{\chi^{-1}(\bar\gamma)}$ for every
$\gamma\in\Gamma$, for which it is sufficient to show that
$\gamma\sigma\gamma^{-1}(x)=x+\chi^{-1}(\bar\gamma)$.

We are given that $\bar\gamma(\bar a)=\chi(\bar\gamma)\bar a$.  As
$\wp(\gamma(x))=\gamma(x)^p-\gamma(x)=\gamma(a)$, we must have
$\gamma(x)=\chi(\bar\gamma)x+b_\gamma$ for some $b_\gamma\in K$, with
$\gamma^{-1}(b_\gamma)+\chi(\bar\gamma)b_{\gamma^{-1}}=0$ (modulo $\wp(K)$)
in order to ensure that $\gamma^{-1}\gamma(x)=x$.  Now, at least modulo
$\wp(K)$, we have
% $$
% \gamma\sigma\gamma^{-1}(x)
% =\gamma\sigma(\chi^{-1}(\bar\gamma)x+b_{\gamma^{-1}})
% =\gamma(\chi^{-1}(\bar\gamma)x+\chi^{-1}(\bar\gamma)+b_{\gamma^{-1}})
% $$
% which equals
% $
% x+\chi^{-1}(\bar\gamma)b_\gamma+\chi^{-1}(\bar\gamma)
%  +\gamma(b_{\gamma^{-1}})
% =x+\chi^{-1}(\bar\gamma)
% $.
$$
\eqalign{
\gamma\sigma\gamma^{-1}(x)
&=\gamma\sigma(\chi^{-1}(\bar\gamma)x+b_{\gamma^{-1}})\cr
&=\gamma(\chi^{-1}(\bar\gamma)x+\chi^{-1}(\bar\gamma)+b_{\gamma^{-1}})\cr
&=x+\chi^{-1}(\bar\gamma)b_\gamma+\chi^{-1}(\bar\gamma)
 +\gamma(b_{\gamma^{-1}})\cr 
&=x+\chi^{-1}(\bar\gamma).
}$$

\bigbreak

We have seen that solvable degree-$p$ extensions $E$ of $F$ give rise to
degree-$p$ cyclic extensions $L$ of $K$ which are galoisian over $F$
(lemma~15) or equivalently to a $G$-stable line $D$ in $K^\times\!/K^{\times
  p}$ or $K^+\!/\wp(K^+)$ (lemma~17), but $E$ is not determined by $D$ or $L$ unless
$E|F$ is cyclic (lemma~12)~; it is determined only up to $F$-conjugacy.

But it should be possible to determine from $L$ or $D$ invariants of $E$ which
depend only on the $F$-conjugacy class of $E$.  We shall see some examples in
the next~\S; here we shall see how to recover the galoisian closure $\tilde E$
of $E$.

Recall (lemma~12) that $\tilde E$ is a degree-$p$ cyclic extension of some
cyclic extension $F'$ of $F$ of degree dividing~$p-1$.  Which degree-$p$
cyclic extensions of $K$ come from some degree-$p$ cyclic extension of a given
$F'$~?

\th LEMMA 18
\enonce
Let\/ $F'$ be an extension of\/ $F$ in $K$ and\/ $G'=\Gal(K|F')$.  A
degree-$p$ cyclic extension $L$ of\/ $K$ of group\/ $P=\Gal(L|K)$ comes from a
degree-$p$ cyclic extension $E'$ of\/ $F'$ if and only if\/ $L$ is galoisian
over\/ $F'$ and the resulting action of\/ $G'$ on\/ $P$ by  conjugation is
trivial. 
\endth
In terms of the line $D$ corresponding to $L$, the condition is that $D$ be
$G'$-stable and that $G'$ should act on $D$ via the cyclotomic character
$\omega$ (lemma~17).  The proof is similar to that of lemma~16 (cf.~{\it
  Remark\/}~4 and \citer\final(),~\S2).

Now let $D$ be in the $\chi$-eigenspace for some character $\chi$ of $G$.  For
$L=K(\!\root p\of D)$ (in characteristic~$\neq p$) or $L=K(\wp^{-1}(D))$ (in
characteristic~$p$) to come from a degree-$p$ cyclic extension of $F'$, a
necessary condition is that $\chi|_{G'}=\omega$ (lemmas~17--18), so the
smallest $F'$ which would do is $F_\chi=K^{G_\chi}$, where
$G_\chi=\Ker(\omega\chi^{-1})$.

Define the short exact sequence $1\to P\to\Gamma_\chi\to G_\chi\to1$ by
restricting $1\to P\to\Gamma\to G\to1$ along the inclusion $G_\chi\to G$.  As
the action of $G_\chi$ on $P$ is trivial (and their orders mutually prime),
there is a canonical section $G_\chi\to\Gamma_\chi$ ({\it Remark\/}~4) using
which we identify $G_\chi$ with a subgroup of $\Gamma$.  Here is the diagram
$$
\def\\{\mskip-2\thickmuskip}
\def\droite#1{\\\hfl{#1}{}{8mm}\\}
\diagram{
1&\rightarrow&P&\droite{}
 &\Gamma_{\phantom{\chi}}&\droite{}&G&\rightarrow&0\phantom.\cr
&&\ufl{=}{}{5mm}&&\ufl{}{}{5mm}&&\ufl{}{\subset}{5mm}\cr
1&\rightarrow&P&\droite{}&\Gamma_\chi&\droite{}&G_\chi&\rightarrow&0.\cr
}
$$
With the identification $G_\chi\subset\Gamma$ we have $EF_\chi=L^{G_\chi}$,
so we have proved~:

\th LEMMA 19
\enonce
Let\/ $E$ be a degree-$p$ solvable extension of $F$ and\/ $\chi$ the
character through which\/ $G$ acts on the corresponding line\/ $D$.  Then\/
$L^{G_\chi}$ is the galoisian closure of\/ $E$, and we have
$L^{G_\chi}=EF_\chi$.
\endth
In summary, we have the following commutative diagram of fields
$$
\def\\{\mskip-2\thickmuskip}
\def\vide{\phantom{phan}}
\diagram{
E&\hfl{I_\chi}{}{8mm}&EF_\chi&\hfl{G_\chi}{}{8mm}&L=EK&\vide
 &G_\chi=\Ker(\omega\chi^{-1})\cr
\ufl{}{}{5mm}&&\ufl{}{P}{5mm}&&\ufl{}{P}{5mm}&\vide
 &P=\Gal(L|K)\cr
F&\hfl{I_\chi}{}{8mm}&F_\chi&\hfl{G_\chi}{}{8mm}&K&\vide
 &I_\chi=\Im(\omega\chi^{-1});\cr
}
$$
in which the names above the arrows are not maps but the relative automorphism
groups (in addition to $\Gamma=\Gal(L|F)$, $\Gamma_\chi=\Gal(L|F_\chi)$ and
$G=\Gal(K|F)$)~; the extension $E|F$ is not galoisian unless $I_\chi=\{1\}$,
which happens precisely when $\chi=\omega$.

\th LEMMA 20
\enonce
The compositum\/ $C$ of all degree-$p$ solvable extensions of\/ $F$ is the
maximal abelian extension $M'$ of exponent dividing\/~$p$ of\/ $K'$, where $K'$
is the compositum of the $F_\chi$ such that the $\chi$-eigenspace
$K^\times\!/K^{\times p}(\chi)$ or\/ $K^+\!/\wp(K^+)(\chi)$ (respectively) is\/
$\neq\{1\}$ or $\neq\{0\}$, where $\chi:G\to\F_p^\times$ runs through
characters of\/ $G=\Gal(K|F)$.
\endth
It is clear that $F_\chi\subset C$ for every $\chi:G\to\F_p^\times$ such that
the $\chi$-eigenspace is $\neq\{1\}$ or $\neq\{0\}$ (lemma~19).  By
definition, $K'$ is the compositum of all such $F_\chi$, so $K'\subset C$.  It
is also clear that that $C\subset M'$ (cf.~lemma~15).  We may thus ask for the
subspace $G_{M'|C}=\Gal(M'|C)$ of the $\F_p$-space $G_{M'|K'}=\Gal(M'|K')$.
Here is the picture of the various fields and their relative automorphism
groups~:
$$
F\overbrace{\underbrace{\longrightarrow}_{I_\chi}F_\chi\longrightarrow}^{G'}K'
\underbrace{\longrightarrow C
  \overbrace{\longrightarrow}^{G_{M'|C}}}_{G_{M'|K'}}M',\qquad
  F\underbrace{\overbrace{\longrightarrow}^{G'}K'\longrightarrow}_G K.
$$
Notice that, as each $F_\chi|F$ is galoisian, so is $K'|F$.  By maximality, so
is $M'|F$.  Now, $C$ is the compositum of all degree-$p$ cyclic extensions $L$
of $K'$ (in $M'$) which are galoisian over $F$~; the set of such $L$ is in
natural bijection $L=M'^H$ with the set of $G'$-stable hyperplanes $H\subset
G_{M'|K'}$, where the action of $G'$ on $G_{M'|K'}$ comes from the short exact
sequence
$$
1\to G_{M'|K'}\to\Gal(M'|F)\to G'\to1,
$$ 
which itself comes from the galoisian tower $M'|K'|F$.  So the subspace in
question $G_{M'|C}\subset G_{M'|K'}$ is the intersection of all $G'$-stable
hyperplanes in $G_{M'|K'}$.  But this intersection is trivial (lemma~8), so
$G_{M'|C}=\{1\}$ and hence $C=M'$.

\smallbreak

{\it Examples 21}\pointir When the field $F$ is finite, we have $K'=F$.
Indeed, the only character $\chi$ of $G$ for which $K^\times\!/K^{\times
  p}(\chi)$ (resp.~$K^+\!/\wp(K^+)(\chi)$) is $\neq\{1\}$ (resp.~$\neq\{0\}$) is
$\chi=\omega$ in characteristic $\neq p$ (resp.~the trivial character $\chi=1$
in characteristic~$p$, because the trace gives an isomorphism $K^+\!/\wp(K^+)=\Fp$),
and we have $F_\chi=F$ in both cases.  When $F$ is a local field of residual
characteristic~$p$, we have $K'=K=F(\!\root{p-1}\of{F^\times})$, and
$C=K(\root p\of{K^\times})$ if $\car(F)=0$ and $C=K(\wp^{-1}(K))$ if
$\car(F)=p$ (prop.~36).  When $F$ is a local field of residual
characteristic~$\neq p$, we have $K'=F(\!\root p\of1)$ and $C=K'(\!\root
p\of{K'^\times})$ (\S8).

\medbreak

{\it Summary of the algebraic ingredients}\pointir We extract from the
foregoing what is relevant for the following \S\S.  We have a prime $p$, a
field $F$, and the maximal abelian extension $K$ of $F$ of exponent dividing
$p-1$~; if the characteristic of $F$ is $\neq p$, then $K$ contains a
primitive $p$-th root of~$1$.  Write $G=\Gal(K|F)$.

Solvable extensions $E$ of $F$ of degree~$p$ give rise to $G$-stable lines $D$
in $K^\times\!/K^{\times p}$ or $K^+\!/\wp(K^+)$, and every $G$-stable line $D$
arises from some $E$.  Two such extensions $E$, $E'$ give rise to the same $D$
if and only if they are $F$-conjugate~; each $E$ has exactly $p$ conjugates,
unless $E|F$ is cyclic.  Invariants of $E$ which depend only on its
$F$-conjugacy class (such as the galoisian closure $\tilde E$) can be
recovered from the $G$-module $D$.  For example, $\tilde E=L^{G_\chi}$, where
$\chi:G\to\F_p^\times$ is the character through which $G$ acts on $D$,
$L=K(\!\root p\of D)$ or $L=K(\wp^{-1}(D))$, and
$G_\chi=\Ker(\omega\chi^{-1})$ has been identified with a subgroup of
$\Gamma=\Gal(L|F)$.  Also, $E|F$ is cyclic if and only if $\chi=\omega$, and
the number of $F$-conjugates of $E$ is $p$ if $\chi\neq\omega$.

\bigbreak 

{\bf 5.  The case of local fields}\pointir We now return definitively to our
local field $F$ with finite residue field $k$ of characteristic $p$ and
cardinality $q=p^f$, and consider degree-$p$ separable extensions of $F$.  (In
\S8, we discuss extensions of $F$ of prime degree $p'\neq p$.)  Let us record
the special features of this case in a series of remarks, of which 23--25
summarise the essential content of \citer\doud().  We denote the normalised
valuation of $F$ by $v:F^\times\to\Z$.

{\it Remark 22}\pointir Every degree-$p$ separable extension $E$ of $F$ is
solvable.  Indeed, the group $\Gamma=\Gamma_{-1}$ of the galoisian closure
$\tilde E|F$ comes equipped with the ramification filtration
$(\Gamma_i)_{i\in\N}$ (in the lower numbering) the successive quotients of
which are commutative.

{\it Remark 23}\pointir It is simpler to show that $\tilde E$ is a degree-$p$
cyclic extension of a cyclic extension $F'$ of $F$ of degree dividing $p-1$~;
we don't need to invoke lemma~2, as lemma~1 suffices.  Indeed, we may assume
that $E|F$ is ramified~; the ramification subgroup $\Gamma_1$ then has order
$>1$, for otherwise $\tilde E|F$ would be tamely ramified whereas the
ramification index of $E|F$ is~$p$.  As the order of ${\goth S}_p$ is not
divisible by $p^2$, $\Gamma_1$ must have order~$p$.  As the subgroup
$\Gamma_1\subset\Gamma$ is normal, the quotient $\Gamma/\Gamma_1$ is a
subgroup of $\F_p^\times$ (lemma~1).
 
{\it Remark 24}\pointir When $E|F$ is ramified, the unique ramification break
$b$ of $\Gamma_1=\Gal(\tilde E|F')$ (the integer $b$ such that
$\Gamma_b=\Gamma_1$, $\Gamma_{b+1}=\{1\}$), and the order $pt$ of the inertia
group $\Gamma_0$ determine the valuation $v(\delta_{E|F})$ of the discriminant
of $E|F$ by computing $v(\delta_{\tilde E|F})$ in two different ways, using
the {\it Schachtelungssatz\/} along the two towers $\tilde E|F'|F$, $\tilde
E|E|F$.  The situation is summarised in the following diagram in which the
ramification indices (resp.~residual degrees) are indicated outside
(resp.~inside) the square
$$
\diagram{
E&\hfl{t}{r}{10mm}&\tilde E\phantom.\cr
\ufl{p}{1}{5mm}&&\ufl{1}{p}{5mm}\cr
F&\hfl{r}{t}{10mm}&F'.\cr
}
$$
Indeed, $v_{F'}(\delta_{\tilde E|F'})=(p-1)(1+b)$.  The extension $F'|F$ is
tamely ramified of group $\Gamma/\Gamma_1$ of some order $tr$ and inertia
subgroup $\Gamma_0/\Gamma_1$ of order~$t$, so the ramification index is $t$
and the residual degree is $r$, whence $v(\delta_{F'|F})=(t-1)r$.  The
extension $\tilde E|E$ has the same ramification index and residual degree as
$F'|F$, so $v_E(\delta_{\tilde E|E})=(t-1)r$.  This leads to the equality
$$
(p-1)(1+b).r+(t-1)r.p=(t-1)r+v(\delta_{E|F}).tr,
$$
which leads to $v(\delta_{E|F})=(p-1)(b+t)/t$.  A similar double application
of the {\it Schachtelungssatz\/} is needed later (prop.~37).

{\it Remark 25}\pointir We are in the presence of {\it two\/} natural
embeddings $\iota,\theta_0:\Gamma_0/\Gamma_1\to\F_p^\times$.  The first one
comes from the conjugation action
$\tau\sigma\tau^{-1}=\sigma^{\iota(\bar\tau)}$ of $\Gamma_0$ on $\Gamma_1$
(and the identification $\F_p^\times=\Aut\Gamma_1$).  The second one comes
from the galoisian action $\bar\tau(\pi)=\theta_0(\bar\tau)\pi$ of
$\Gamma_0/\Gamma_1$ on the $t$-th roots $\pi$ of a uniformiser $\varpi$ of the
maximal unramified extension $F'_0$ of $F$ in $F'$ such that $\varpi\in
F'^{\times t}$ (and the identification $\F_p^\times\subset F_0'^\times$ with
the group of $(p-1)$-th roots of~1)~; $\theta_0$ is independent of the choice
of $\varpi$.  We also have an embedding $\theta_b:\Gamma_b/\Gamma_{b+1}\to
U_b/U_{b+1}$ (where $U_i$ ($i>0$) is the group of principal units of $F'$ of
level at least $i$).  But $\Gamma_b=\Gamma_1$ and $\Gamma_{b+1}=\{1\}$, so
$\theta_b$ is an embedding of $\Gamma_1$.  Now, the compatibility relation
$\theta_b(\tau\sigma\tau^{-1}) =\theta_b(\sigma)^{\theta_0(\bar\tau)^b}$ for
every $\sigma\in\Gamma_1$ and every $\tau\in\Gamma_0$ \citer\corpslocaux(IV,
prop.~9) implies that $\iota(\bar\tau)=\theta_0(\bar\tau)^b$ for every
$\bar\tau\in\Gamma_0/\Gamma_1$.  In other words, the two embeddings
$\iota,\theta_0$ differ by the automorphism $(\ )^b$ of $\Gamma_0/\Gamma_1$.
In particular, $\gcd(b,t)=1$ (where $t$ is the order of $\Gamma_0/\Gamma_1$).

{\it Remark 26}\pointir In principle, it should now be possible to prove
Serre's degree-$p$ mass formula by computing the contribution of each such
$F'$~; when $F'=F$, then $E|F$ is cyclic, and the contribution of these has
been computed in \citer\final(prop.~14--16).  The number of $F'$ can be
deduced from \citer\hasse(Kap.~16) or \citer\course(Lecture~18), and equals
the number of cyclic subgroups of $F^\times\!/F^{\times p-1}$.

{\it Remark 27}\pointir The maximal abelian extension $K$ of $F$ of exponent
dividing $p-1$ equals $F(\!\root p-1\of{F^\times})$.  Indeed, $F$ contains a
primitive $(p-1)$-th root of~$1$.  The group $G=\Gal(K|F)$ is dual to
$F^\times\!/F^{\times p-1}$ under the pairing
$$
G\times (F^\times\!/F^{\times  p-1})\to\F_p^\times
\quad(\sigma,\bar x)\mapsto{\sigma(y)\over y}\ \ (y^{p-1}=x)
$$
in which $\F_p^\times\subset F^\times$ has been identified with the group
of $(p-1)$-th roots of~$1$.  

In the characteristic-$0$ case, denote by $\omega:G\to\F_p^\times$ the
cyclotomic character giving the action of $G$ on the $p$-th roots of~$1$~; in
characteristic~$p$, let $\omega=1$ be the trivial character.  Note that $K$ is
a finite extension of $F$ of ramification index and residual degree $p-1$.
Denote its ring of integers by $\ogoth$, the unique maximal ideal of $\ogoth$
by $\pgoth$, and the residue field by $l=\ogoth/\pgoth$.  Finally, let
$U_n=1+\pgoth^n$ for $n>0$.

{\it Remark 28}\pointir In characteristic $0$, the character $\omega$
corresponds to the class $\overline{-p}\in F^\times\!/F^{\times p-1}$ under
biduality.  In other words, we have to show that
$\sigma(y)/y\equiv\omega(\sigma)\pmod{\pgoth}$ for every $\sigma\in G$ and
every $(p-1)$-th root $y\in K^\times$ of~$-p$.  Let $\zeta\in K^\times$ be a
primitive $p$-th root of~$1$, so that $\sigma(\zeta)=\zeta^{\omega(\sigma)}$
for every $\sigma\in G$.  We may take $y=\eta\pi$, where $\pi=1-\zeta$ and
$\eta\in U_1$~; we have $\sigma(\pi)/\pi\equiv\omega(\sigma)\pmod{\pgoth}$,
from which the claim follows.  Cf.~\citer\locdisc(prop.~25).

Remark finally that the space $K^\times\!/K^{\times p}$ or $K^+\!/\wp(K^+)$ carries
a natural filtration, and that the discriminant of a degree-$p$ separable
extesnion $E$ of $F$ can be computed from the ``\thinspace level\thinspace''
of the corresponding line $D\subset K^\times\!/K^{\times p}$ (resp.~$D\subset
K^+\!/\wp(K^+)$) in this filtration.  See the next~\S\ for a description of the
filtration and the definition of the ``\thinspace level\thinspace'' of a line,
and prop.~37 for the computation.

\bigbreak

{\bf 6.  Filtered galoisian modules}\pointir We keep the notation $F,K, G$
from~\S5. We have seen (\S4) that the $G$-modules $K^\times\!/K^{\times p}$ or
$K^+\!/\wp(K^+)$ (respectively) play an important role in the study of degree-$p$
separable extensions of $F$.  These modules come with a natural filtration,
which we discuss next.

Denote by $(\bar U_n)_{n>0}$ the filtration on $\bar U_0=K^\times\!/K^{\times
  p}$ by units of various levels~.  Similarly, in the characteristic-$p$ case,
let ~$\overline{{\goth p}^n}$ be the image of $\pgoth^n$ in the
$\Fp$-space~$K^+\!/\wp(K^+)$, where $\pgoth$ is the unique maximal ideal of the ring
of integers~$\ogoth$ of $K$.  For some background on these $\F_p$-space
(without the $G$-action), see \citer\locdisc() in characteristic~$0$ and
\citer\further() in characteristic~$p$.

Examples of $G$-stable lines are provided, in the characteristic-$0$ case, by
$\bar U_{pe}$ and $\bar\mu$ (the image of the torsion subgroup
$\mu\subset{\goth o}^\times$), on both of which $G$ acts via the cyclotomic
character $\omega$.  In the characteristic-$p$ case, the line $\bar{\goth
  o}=\Fp$ is $G$-stable and the action of $G$ is in fact trivial.

In general, the subspaces $\bar U_i$ for $i\in[0,pe]$ in characteristic~$0$
(resp.~$\overline{{\goth p}^i}$ for $i\in-\N$ in characteristic~$p$) are
$G$-stable, essentially because there is a unique extension of the valuation
from $F$ to $K$.  We have seen that $\bar U_{pi+1}=\bar U_{pi}$ except for
$i=0,e$ (\citer\locdisc(prop.~42)) and that
$\overline{\pgoth^{pi+1}}=\overline{\pgoth^{pi}}$ except for $i=0$
(\citer\further(prop.~11)), respectively.

The codimension is~$1$ in the three exceptional cases, namely $\bar
U_1\subset\bar U_0$ and $\bar U_{pe_1+1}\subset\bar U_{pe_1}$ in
characteristic~$0$, and $\overline{\pgoth}\subset\overline{\ogoth}$ in
characteristic~$p$.  In characteristic~$0$, we have $\bar U_{pe+1}=\{1\}$ and
$\bar U_{pe}$ is a stable line on which $G$ acts via $\omega$, and the
valuation $v_K$ provides an isomorphism $\bar U_0/\bar U_1\to\Z/p\Z$.  In
characteristic~$p$, we have $\bar\pgoth=0$ and the trace $S_{l|\Fp}$ induces
an isomorphism $\bar{\goth o}\to\Fp$, where $l$ is the residue field of $K$.

For $i\in[1,pe[$ prime to~$p$ (in characteristic~$0$) or $i<0$ prime to~$p$
(in characteristic~$p$), the codimension equals the absolute degree $[l:\Fp]$,
and indeed the quotients are canonically isomorphic to $U_i/U_{i+1}$
(resp.~${\goth p}^{i}/{\goth p}^{i+1}$) (\citer\locdisc(prop.~42),
\citer\further(prop.~11)).  Thus they are not merely $\Fp$-spaces but
$k$-spaces (of $k$-dimension $p-1$).  The pictures in \citer\final(), \S5
summarise some of these facts.

\th PROPOSITION 29
\enonce
For\/ $i\in[1,pe[$ prime to\/~$p$ (in characteristic\/~$0$) or for\/ $i<0$
prime to\/~$p$ (in characteristic\/~$p$), the natural maps\/ 
$$
\bar U_i/\bar U_{i+1}\to U_i/U_{i+1}\to\pgoth^i/\pgoth^{i+1},\quad
(\hbox{resp. }\overline{\pgoth^i}/\overline{\pgoth^{i+1}} \to
\pgoth^i/\pgoth^{i+1})
$$
of\/ $k$-spaces are $G$-equivariant $k$-isomorphisms.
\endth 
{\it Beweis~:} Klar.

% Indeed, the isomorphisms in the first assertion are induced by the identity
% map on $U_i$ (resp.~${\goth p}^i$)~: the diagrams
% $$
% \def\\{\mskip-2\thickmuskip}
% \def\droite#1{\\\hfl{#1}{}{8mm}\\}
% \def\vide{\phantom{phantom}}
% \diagram{
% U_i&\droite{=}&U_i&\vide&\pgoth^i&\droite{=}&\pgoth^i\cr
% \vfl{}{}{5mm}&&\vfl{}{}{5mm}&\vide&\vfl{}{}{5mm}&&\vfl{}{}{5mm}\cr
% \bar U_i/\bar U_{i+1}&\droite{\sim}&U_i/U_{i+1}
%  &\vide&\overline{\pgoth^i}/\overline{\pgoth^{i+1}}&\droite{\sim}
%  &\pgoth^i/\pgoth^{i+1}\cr
% }
% $$ 
% are commutative.  The second assertion is also clear, as the isomorphism in
% question is nothing but $\bar u\mapsto\overline{u-1}$.

\smallskip

So we need to study the $G$-modules $\pgoth^i/\pgoth^{i+1}$, for which a
preliminary study of the $G$-module $l$ is useful.

\th PROPOSITION 30
\enonce
Let\/ $k$ be a finite field, $q=\Card k$, and\/ $l|k$ any extension of degree 
dividing $q-1$.  For every character\/ $\chi:\Gal(l|k)\to k^\times$, the
$\chi$-eigenspace $l(\chi)$ is a $k$-line in $l$.
\endth
Although this is a special case of lemma~9, we give a short direct proof.  Let
$\varphi$ (Frobenius) be the canonical generator $x\mapsto x^q$ of
$\Gal(l|k)$, and put $a=\chi(\varphi)$~; let $m=[l:k]$ be the order of
$\varphi$, so that the order of $a\in k^\times$ divides~$m$.  The
$\chi$-eigenspace consists of all $x\in l$ such that $\varphi(x)=ax$~; such
$x$ are roots of the binomial $T^q-aT=T(T^{q-1}-a)$, which has at most $q$
roots.  As the map $x\mapsto x^q-ax$ defined by this binomial is a linear
endomorphism of the $k$-space $l$, it is sufficient to prove that $a$ has a
$(q-1)$-th root in $l$.

We have said that the order of $a$ in $k^\times\!/k^{\times q-1}=k^\times$
divides $m$. Therefore the degree of the extension $k(\!\root{q-1}\of a)|k$
divides $m$, and hence $k(\!\root{q-1}\of a)\subset l$.  This shows that the
$k$-endomorphism $x\mapsto x^q-ax$ of $l$ is not injective, and hence its
kernel $l(\chi)$ is a $k$-line in $l$.  Incidentally, if $\chi$ is trivial,
then $a=1$ and $l(\chi)=k$.
 
\smallskip

Momentarily let $K$ be any galoisian extension of $F$ of group $G$, and
suppose that there is a uniformiser $\pi$ of $K$ such that $\varpi=\pi^s$ is
in $F$ for some $s>0$.

\th PROPOSITION 31
\enonce
For every integer\/ $i\in\Z$, ``\thinspace multiplication by\/
$\varpi$'' gives an isomorphism\/
$\pgoth^i\!/\pgoth^{i+1}\to\pgoth^{i+s}\!/\pgoth^{i+s+1}$ of\/ $k[G]$-modules.
\endth
More precisely, the reduction modulo $\pgoth$ of the $\ogoth$-linear
isomorphism $x\mapsto\varpi x:\pgoth^i\to\pgoth^{i+s}$ is $G$-equivariant.
But this is clearly the case~: $\sigma(\varpi x)=\varpi\sigma(x)$ for every 
$\sigma\in G$, because $\varpi\in F$ and $\sigma$ is $F$-linear.

\smallskip

Let us now return to our $K=F(\!\root{p-1}\of{F^\times})$ and
$G=\Hom(F^\times\!/F^{\times p-1},\F_p^\times)$, so that the group of
characters of $G$ is $\Hom(G,\F_p^\times)=F^\times\!/F^{\times p-1}$.  (The
cyclotomic character $\omega$ corresponds to $\overline{-p}$ in
characteristic~$0$, and $\omega=1$ by convention in characteristic~$p$).  Each
character $\chi$ therefore has a ``\thinspace valuation\thinspace'' $\bar
v(\chi)\in\Z/(p-1)\Z$, coming from the valuation $v:F^\times\to\Z$.

Unramified characters --- those which are trivial on the inertia subgroup
$G_0$ --- are the same as characters of the quotient $G/G_0=\Gal(l|k)$.  They
get identified with the kernel $\ogoth_F^\times\!/\ogoth_F^{\times
  p-1}=k^\times\!/k^{\times p-1}$ of $\bar v$.  Indeed, the short exact
  sequence $1\to G_0\to G\to G/G_0\to1$ gives rise, upon taking duals $(\ 
  )^\vee=\Hom(\ ,\F_p^\times)$, to a short exact sequence which fits into the
  commutative diagram
$$
\def\\{\mskip-2\thickmuskip}
\def\droite#1{\\\hfl{#1}{}{8mm}\\}
\diagram{
0&\rightarrow&(G/G_0)^\vee&\droite{}&
 G^\vee&\droite{}&G_0^\vee&\rightarrow&0\phantom.\cr
&&\vfl{\sim}{}{5mm}&&\vfl{\sim}{}{5mm}&&\vfl{\sim}{}{5mm}\cr
1&\rightarrow&k^\times\!/k^{\times p-1}&\droite{}
 &F^\times\!/F^{\times p-1}&\droite{\bar v}&\Z/(p-1)\Z&\rightarrow&0.\cr
}
$$
Notice that the map $\bar x\mapsto x^{(q-1)/(p-1)}$ ($x\in k^\times$)
identifies $k^\times\!/k^{\times p-1}$ with $\F_p^\times$.  In
characteristic~$0$, the cyclotomic character $\omega$ is unramified if and
only if $p-1|e$, for $\omega$ corresponds to $\overline{-p}$ and $\bar
v(\overline{-p})\equiv e\pmod{p-1}$.

Let us decompose the $k[G]$-modules $\pgoth^i/\pgoth^{i+1}$ as an
internal direct sum of $\chi$-eigenspaces for various
$\chi:G\to\F_p^\times$.  For $i=0$, the $k[G]$-module
$\pgoth^i/\pgoth^{i+1}$ is in fact the $k[G/G_0]$-module $l$.  We have
seen that for every unramified character $\chi$ of $G$, the
$\chi$-eigenspace $l(\chi)$ is a $k$-line (lemma~9 or prop.~30).  It
follows that for every {\it ramified\/} character $\chi$, we have
$l(\chi)=0$~: there is room only for so many, and unramified
characters have used it all up.

Let us provide the details of the notion of {\it twisting\/} a $k[G]$-module
$m$ by a character $\xi:G\to\F_p^\times$.  Denote by $m\{\xi\}$ the
$k[G]$-module whose underlying $k$-space is $m$, but the new action
$\star_\xi$ is defined by $\sigma\star_\xi x=\xi(\sigma)\sigma(x)$ for every
$\sigma\in G$ and every $x\in m$, so that if $\xi=1$ is the trivial character,
then $\sigma\star_1x=\sigma(x)$ and $m\{1\}$ is the $k[G]$-module $m$.  It is
clear that $m\{\xi_1\xi_2\}=m\{\xi_1\}\{\xi_2\}$ for any two characters
$\xi_1$, $\xi_2$ of $G$, and that $(m_1\oplus m_2)\{\xi\}=m_1\{\xi\}\oplus
m_2\{\xi\}$ for any two $k[G]$-modules $m_1$, $m_2$.  In this process, the
$\chi$-eigenspace of $m$ gets converted into the $\xi\chi$-eigenspace of
$m\{\xi\}$, for every character $\chi$ of $G$.  The same discussion
applies to $\F_p[G]$-modules.

It is easy to see that $l\{\xi\}$ is $k[G]$-isomorphic to $l$ for every {\it
  unramified\/} character $\xi$ of $G$.  Indeed, by prop.~30, $l$ is the
direct sum of the $k$-lines $l(\chi)$ indexed by the unramified characters
$\chi\in\Hom(G/G_0,\F_p^\times)$ of $G$, so $l\{\xi\}$ is the direct sum of
the $k$-lines $l(\xi\chi)$.  But as $\chi$ runs through
$\Hom(G/G_0,\F_p^\times)$, so does $\xi\chi$, and the two direct sums are
$k[G]$-isomorphic.

It follows that $l\{\xi\}$ depends only on $\bar v(\xi)\in\Z/(p-1)\Z$ (up to
$k[G]$-isomorphism).  We denote by $l[i]$ (for $i\in\Z$) the $k[G]$-modules
$l\{\xi\}$ for any $\xi$ such that $\bar v(\xi)\equiv i\pmod{p-1}$, so that
$l[0]=l$.

\th PROPOSITION 32
\enonce
For every\/ $i\in\Z$ and every character\/ $\chi:G\to\F_p^\times$, the\/
$\chi$-eigenspace in the\/ $k[G]$-module\/ $l[i]$ is a\/ $k$-line if\/ $\bar
v(\chi)\equiv i\pmod{p-1}$~; it is reduced to\/~$0$ otherwise.
\endth
This is just prop.~30 in the case~$i=0$, and the general case follows from
this by our discussion of twisting.  An immediate consequence is the following
result. 

\th PROPOSITION 33
\enonce
For\/ $i\in\Z$, the\/ $k[G]$-module\/ $\pgoth^i\!/\pgoth^{i+1}$ is
isomorphic to\/ $l[i]$.
\endth
Choose any uniformiser $\pi$ of $K$ such that $\pi^{p-1}$ is (a uniformiser)
in $F$~; this is possible.  It is easy to see that, by taking $\pi^i$ as an
$\ogoth$-basis of $\pgoth^i$, the resulting $k$-linear map
$\pgoth^i/\pgoth^{i+1}\to l\{\xi^i\}$ is an isomorphism of $k[G]$-modules,
where $\xi$ is the character such that $\sigma(\pi)=\xi(\sigma)\pi$ for every
$\sigma\in G$.  As we have $\bar v(\xi)\equiv1$, this
shows that $\pgoth^i/\pgoth^{i+1}$ is $k[G]$-isomorphic to $l[i]$.

\th COROLLARY 34
\enonce
Put\/
$W_1=\pgoth^{-1}\!/\pgoth^{0}\oplus
\pgoth^{-2}\!/\pgoth^{-1}\oplus\cdots 
\oplus\pgoth^{-(p-1)}\!/\pgoth^{-(p-1)+1}$. 
For every\/ $\chi:G\to\F_p^\times$, the\/ $\chi$-eigenspace in
the\/ $k[G]$-module\/ $W_1$ is a\/ $k$-line. 
\endth
In other words, and thanks to prop.~33, $W_1=l[-1]\oplus\cdots\oplus
l[-(p-1)]$.  The point is that we have endowed $W_1$ with the filtration for
which the successive quotients are, for example when $p=5$,
$$
l[-1],\ l[-2],\ l[-3],\ l[-4],
$$
in that specific order, rather than in some other order.  If we twist it by
the cyclotomic character (which we will soon need to) to get $W_1\{\omega\}$,
and if $\omega$ has ``\thinspace valuation\thinspace'' $\bar
v(\omega)\equiv1\pmod4$ (for example when $F|\Q_5$ is unramified), then the
successive quotients of the filtered $k[G]$-module $W_1\{\omega\}$ are
$$
l[-4],\ l[-1],\ l[-2],\ l[-3],
$$
which are circularly shifted one step to the right (or three steps to the
left, which comes to the same because $1\equiv-3\pmod4$).  We need to keep
track of both the filtration and the $G$-action.  There is no difference when
$\omega$ is unramified, for then the shift is by~$0$ steps.

In a similar vein, define $W_2=\pgoth^{-(p+1)}\!/\pgoth^{-p}\oplus \cdots
\oplus\pgoth^{-(2p-1)}\!/\pgoth^{-(2p-1)+1}$, and think of this $k[G]$-module
as being endowed with the filtration for which the successive quotients are
$l[-(p+1)],\ldots,l[-(2p-1)]$ (cf.~prop.~33), in that specific order, so that
when $p=5$, they are
$$
l[-2],\ l[-3],\ l[-4],\ l[-1].
$$
As in prop.~34, we see that the\/ $\chi$-eigenspace in the\/ $k[G]$-module\/
$W_2$ is a\/ $k$-line for every\/ $\chi:G\to\F_p^\times$.

For every $i\in\N$, put $W_{i+1}=l[-(ip+1)]\oplus\cdots\oplus l[-(ip+p-1)]$.
These $k[G]$-modules are all isomorphic to each other (and free of rank~$1$,
to boot) because the $\chi$-eigenspace in each $W_{i+1}$ is a\/ $k$-line for
every character $\chi:G\to\F_p^\times$, but they carry different filtrations.
The {\it filtered\/} $k[G]$-modules $W_m$, $W_n$ are isomorphic if and only if
$m\equiv n\pmod{p-1}$.

Get back to the $k[G]$-modules $\bar U_i/\bar U_{i+1}$
(resp.~$\overline{\pgoth^i}/\overline{\pgoth^{i+1}}$) for appropriate~$i$.

\th PROPOSITION 35
\enonce
The\/ $k[G]$-module\/ $\bar U_i/\bar U_{i+1}$ for\/ $0<i<pe$ prime to\/~$p$ in
the characteristic-$0$ case
(resp.~$\overline{\pgoth^i}/\overline{\pgoth^{i+1}}$ for\/ $i<0$ prime
to\/~$p$ in the characteristic-$p$ case) is isomorphic to\/  $l[i]$.
\endth
Indeed, $\bar U_i/\bar U_{i+1}$
(resp.~$\overline{\pgoth^i}/\overline{\pgoth^{i+1}}$) is $k[G]$-isomorphic to
$\pgoth^i\!/\pgoth^{i+1}$ (prop.~29), which is $k[G]$-isomorphic to $l[i]$
(prop.~33). 

\smallskip

Let us record all this information in pictures, for I have still not got over
the fact that the quotients are $k[G]$-modules for appropriate~$i$, instead of
merely being $\Fp[G$]-modules.  In characteristic~$0$, the first picture of
the filtered $\Fp[G]$-module $K^\times\!/K^{\times p}=\bar U_0$ is
$$\eqalign{\{\bar1\}\underbrace{\subset}_{\Fp\{\omega\}}\bar
U_{pe}&\underbrace{\subset}_{l[pe-1]}\bar U_{pe-1}\cdots\cr
&\cdots\underbrace{\subset}_{l[pj+1]}\bar U_{pj+1}
\underbrace{=}_{\{1\}}\bar U_{pj}\underbrace{\subset}_{l[pj-1]}\bar U_{pj-1}\cdots\cr
&\phantom{\cdots\underbrace{\subset}_{l[pj+1]}
 =\bar U_{pj}\underbrace{\subset}_{l[pj-1]}\bar U_{pj-1}}
\cdots\bar U_2\underbrace{\subset}_{l[1]}\bar U_1
\underbrace{\subset}_{\Fp}\bar U_{0},\cr
}
$$
with successive quotients indicated below the inclusion signs and
$j\in[1,e[$, whereas in characteristic~$p$ the picture of the $\Fp[G]$-module
$K^+\!/\wp(K^+)$ goes on for ever ($j<0$)~:
$$
\{\bar0\}\underbrace{\subset}_{\Fp}
\bar\ogoth\underbrace{\subset}_{l[-1]}
\overline{\pgoth^{-1}}
\cdots
\underbrace{\subset}_{l[pj+1]}\overline{\pgoth^{pj+1}}
\underbrace{=}_{\{0\}}\overline{\pgoth^{pj}}
\underbrace{\subset}_{l[pj-1]}
\overline{\pgoth^{pj-1}}
\cdots\subset K^+\!/\wp(K^+)
$$
The beauty of this can reduce even the most hardened criminal to tears.

The analogy can be further improved.  First, declare $\omega$ to be the
trivial character in the characteristic-$p$ case.  Secondly, the two pictures
will look even more similar if the first one is shifted to the right by $pe$
steps.  The problem is that the $k[G]$-modules $l[pe-1]$ and $l[-1]$ are not
isomorphic, unless $\omega$ is unramified, which is equivalent to
$e\equiv0\pmod{p-1}$  More precisely, $l[pe-1]$ is the sum of $G$-stable
$k$-lines indexed by the characters of ``\thinspace valuation\thinspace''
$\bar v(\omega)-1$, whereas $l[-1]$ is the sum of $G$-stable $k$-lines indexed
by the characters of ``\thinspace valuation\thinspace'' $-1$.  But this can be
easily remedied if we twist the latter by $\omega$.

So $l[-1]\{\omega\}$ is the same $k[G]$-module as $l[pe-1]$.  Suppressing the
terms indexed by multiples of~$p$ other than $0$ (and $pe$ in
characteristic~$0$), and exploiting the fact that any $G$-stable subspace has
a $G$-stable supplement (lemma~8), the filtered $\Fp[G]$-module
$K^\times\!/K^{\times p}$ (resp.~$K^+\!/\wp(K^+)$) is
$$
\Fp\{\omega\}\oplus l[-1]\{\omega\}\oplus\cdots 
\oplus l[-b^{(n)}]\{\omega\}\oplus\cdots\ \big(\oplus\Fp\big)
$$
where the middle terms are indexed by the sequence $b^{(n)}$ of
\hbox{prime-to-$p$} integers, for every $n>0$ in characteristic~$p$ but only
for $n\in[1,(p-1)e]$ in characteristic~0, which is also when the last
parenthetical term appears.  This picture keeps track of both the filtration
and the $G$-action.

Recall that $b^{(n)}=n+\lfloor(n-1)/(p-1)\rfloor$, and that if $n=(p-1)i+j$
with $i\in\N$ and $j\in[1,p-1]$ ({\it sic\/}), then $b^{(n)}=pi+j$.  Clearly,
$b^{(\ )}:\N^*\to\N^*$ is the unique strictly increasing function whose image
is the set of integers $>0$ prime to~$p$~; we put $b^{(0)}=0$.

{\it Define\/} $V_i$ ($i>0$) to be the filtered $k[G]$-module $W_i\{\omega\}$,
so that for example when $p=5$ and $e=1$, the successive quotients of $V_2$
are
$$
l[-1],\ l[-2],\ l[-3],\ l[-4].
$$

Grouping together $p-1$ middle terms at a time, and recalling the definitions
of $W_i$ and $V_i=W_i\{\omega\}$, we see from the above description that the
filtered $\Fp[G]$-module $K^\times\!/K^{\times p}$ is
$\Fp[G]$-isomorphic to~\footnote{(*)}{(2011/07/14) In the earlier versions, it
  was asserted that the spaces in question are isomorphic to
  $\Fp\{\omega\}\oplus V_1^e\oplus\Fp$ (resp.~$\Fp\oplus V_1\oplus V_1\oplus
  V_1\oplus\cdots$) as filtered $\Fp[G]$-modules.  This was wrong but did not
  affect the proof in \S7, for the lapse merely had the effect of permuting
  the terms in the sum to be computed.}
$$
\Fp\{\omega\}\oplus (V_1\oplus0)\oplus\cdots\oplus(V_{e-1}\oplus0)\oplus V_e\oplus\Fp
$$
if $e<+\infty$ (resp.~$K^+\!/\wp(K^+)$ is isomorphic to
$\Fp\oplus(V_1\oplus0)\oplus(V_2\oplus0)\cdots$ if $e=+\infty$), where
we have inserted some $0$s to indicate the cases of equality.  From
now on, such subtilities will be ignored.

Henceforth we use the {\it opposite\/} filtration on these $\Fp[G]$-modules
\footnote{(*)}{(2011/07/04) This convention was not adopted in the version of
  this paper publised in the {\it Monatshefte}, hence the minus sign in
  prop.~35 there and the absence of this sign in prop.~37 below.}, so that it
becomes an increasing filtration indexed by $[0,pe]$ in characteristic~$0$ and
by $\N$ in characteristic~$p$.  Concretely, the filtration
$$\F_p\{\omega\}=
{\cal   F}_0\subset
{\cal F}_1\subset\cdots
\subset{\cal F}_{pe-1}\subset
{\cal F}_{pe}=K^\times\!/K^{\times p}
$$
in characteristic~$0$ is just ${\cal F}_i=\bar U_{pe-i}$, which justifies
the convention $\bar U_0=K^\times\!/K^{\times p}$.  In characteristic~$p$, the
filtration ${\cal F}_0\subset{\cal F}_1\subset\cdots$ on $K^+\!/\wp(K^+)$ is just
${\cal F}_i=\overline{\pgoth^{-i}}$.  In both cases, ${\cal F}_{pt}={\cal
  F}_{pt-1}$ for every $t\in[1,e[$ and indeed for every $t>0$ in
characteristic~$p$.  Also, ${\cal F}_{pt}=\Fp\{\omega\}\oplus V_1\oplus\cdots
\oplus V_t$ for every $t\in[1,e[$ (resp.~$t>0$).

Such isomorphisms hold because for example $\bar U_{pe-1}/\bar U_{pe}$ is
isomorphic to $l[pe-1]$ (prop.~29 and~33), which is isomorphic to
$l[e-1]=l[-1]\{\omega\}$ (recalling that $\bar v(\omega)\equiv e\pmod{p-1}$),
so ${\cal F}_p/{\cal F}_0=W_1\{\omega\}=V_1$ as filtered $k[G]$-modules.
Also, when $e>1$, 
$$
{\cal F}_{p+1}/{\cal F}_p=
\cases{\pgoth^{pe-(p+1)}/\pgoth^{pe-p}\cr
 \pgoth^{-(p+1)}/\pgoth^{-p}\cr}
=\cases{l[e-2]& if $\car(F)=0$,\cr\;\; l[-2]& if $\car(F)=p$,\cr}
$$
so ${\cal F}_{2p}=\Fp\{\omega\}\oplus W_1\{\omega\}\oplus
W_2\{\omega\}=\Fp\{\omega\}\oplus V_1\oplus V_2$, and so on. Maybe we could
have done with fewer letters, but there we are.

Define the {\it level\/} of a line $D\subset K^\times\!/K^{\times p}$ or
$D\subset K^+\!/\wp(K^+)$ to be the integer $i$ such that $D\subset{\cal F}_i$ but
$D\not\subset{\cal F}_{i-1}$.  The possible levels in characteristic~0 are
$b^{(n)}$ for every $n\in[0,(p-1)e]$, and $pe$.  In characteristic~$p$, they
are $b^{(n)}$ for every $n\in\N$.

\th PROPOSITION 36
\enonce
The compositum\/ $C$ of all degree-$p$ separable extensions of\/ $F$ is the
maximal abelian extension $M$ of exponent\/~$p$ of\/ $K=F(\!\root{p-1}\of
{F^\times})$, namely\/ $M=K(\!\root p\of{K^\times})$ or\/ $M=K(\wp^{-1}(K))$
respectively. 
\endth
This is the main result of \citer\delcorso() in the characteristic-$0$ case,
and their proof can now be carried over to \hbox{characteristic~$p$}.

In view of lemma~20, all we need to show is that $K$ is the compositum of the
$F_\chi$ such that the $\chi$-eigenspace $K^\times\!/K^{\times p}(\chi)$ or\/
$K^+\!/\wp(K^+)(\chi)$ (respectively) is\/ $\neq\{1\}$ or $\neq\{0\}$.  The above
description shows that such is the case for every $\chi:G\to\F_p^\times$, and
it is clear that $K$ is the compositum of all cylic extensions of $F$ of
degree dividing $p-1$.

\smallskip

Interesting as it is, this result will not be needed in the next \S~; only the
algebraic ingredients summarised at the end of \S4 and the structural analysis
of this~\S\ will be needed. \footnote{(*)}{(2011/07/03) Notice that in
  characteristic~$0$, the degree $[M:F]$ is finite and hence $F$ has only
  finitely many extensions of degree~$p$, which implies that it has only
  finitely many extensions of any given degree.}

\bigbreak

{\bf 7. Serre's mass formula in prime degree}\pointir It is time to do the
counting.  The reader will need to refer back frequently to the pictures in
\S6 in what follows.

Let $E|F$ be a ramified degree-$p$ separable extension, $L=EK$, and let
$$
D\subset\cases{\Fp\{\omega\}\oplus V_1\oplus\cdots\oplus V_e\oplus\Fp
 &if $\car(F)=0$,\cr
 \Fp\oplus V_1\oplus V_2\oplus V_3\oplus\cdots
 &if $\car(F)=p$,\cr}
$$
be the $G$-stable line corresponding to $L|K$, so that $D$ is the line we
associate to $E$.  Denote the ``\thinspace level''\thinspace of $D$ by $d(D)$,
so that $d(D)=b^{(n)}$ for some $n\in[1,(p-1)e]$ (resp.~$n>0$) in the {\it peu
  ramifi{\'e}\/} case and $d(D)=pe$ in the {\it tr{\`e}s ramifi{\'e}\/} case
(which occurs only in characteristic~$0$).  Denote by $b(L|K)$ the unique
ramification break of $L|K$, and recall that $c(E)=v(\delta_{E|F})-(p-1)$
measures the wild ramification of $E|F$.

\th PROPOSITION 37
\enonce
We have\/ $b(L|K)=c(E)=d(D)$.
\endth
We have already seen that $b(L|K)=d(D)$ in \citer\locdisc(prop.~60) in the
characteristic-$0$ case (note that the numbering convention here is shifted by
$pe$ steps and involves a sign change) and in \citer\further(prop.~14) in the
characteristic-$p$ case.  Using the {\it Schachtelungssatz} twice in the
diagram
$$
\diagram{
E&\hfl{p-1}{p-1}{10mm}&L\cr
\ufl{p}{1}{5mm}&&\ufl{1}{p}{5mm}\cr
F&\hfl{p-1}{p-1}{10mm}&K\cr
}
$$
in which the ramification indices (resp.~residual degrees) are indicated
outside (resp.~inside) the square, we get $c(E)=b(L|K)$. Cf.~{\it
  Remark\/}~24. 

\th COROLLARY 38
\enonce
The invariant\/ $c(E)$ is prime to $p$ except for\/ $c(E)=pe$ in
characteristic~$0$.
\endth
Indeed, $d(D)$ is prime to $p$ for every line $D\subset K^+\!/\wp(K^+)$ other than
$\bar\ogoth$ \citer\further(prop.~11), and also for every $D\subset
K^\times\!/K^{\times p}$ other than $\bar U_{pe}$, unless $D\not\subset\bar
U_1$, in which case $d(D)=pe$ \citer\locdisc(prop.~42).  It is also clear that
every $m>0$ (resp.~$0<m<pe$) prime to~$p$ occurs as $c(E)$ for some (ramified)
$E$. 

\bigbreak

{\it Proof of the mass formula in degree~$p$}\pointir The number of $E$ which
give rise to the same $D$ is~$1$ if the character $\chi$ through which $G$
acts on $D$ is $\omega$, and $p$ if $\chi\neq\omega$ (\S4).  So the
contribution of such $E$ to Serre's mass formula is (prop.~37)
$$
\sum_{E\mapsto D}q^{-c(E)}=\cases{
\phantom{p}q^{-d(D)}&if $\chi=\omega$,\cr
pq^{-d(D)}&if $\chi\neq\omega$.\cr}
$$
Thus the sum over all ramified separable degree-$p$ extensions $E$ of $F$ gets
replaced by a sum over all $G$-stable lines $D$ in $K^\times\!/K^{\times p}$
or $K^+\!/\wp(K^+)$ other than the level-$0$ line $\bar U_{pe}=\F_p\{\omega\}$
(resp.~$\bar\ogoth=\F_p$) in characteristic~$0$ (resp.~$p$).

For every character $\chi:G\to\F_p^\times$, the dimension of the
$\chi$-eigenspace $(\Fp\{\omega\}\oplus V_1\oplus\cdots\oplus V_i)(\chi)$ (for
all $i\in\N$ in characteristic~$p$, for $i\in[0,e]$ in characteristic~$0$) is
(\S6)
$$
\dim_{\Fp}(\Fp\{\omega\}\oplus V_1\oplus\cdots\oplus V_i)(\chi)
=\cases{
1+if&if $\chi=\omega$,\cr
\phantom{1+\ }if&if $\chi\neq\omega$,\cr}
$$
so the number of points in this space is $pq^i$ if $\chi=\omega$ and $q^i$
if $\chi\neq\omega$.  Therefore the number of lines in $(\Fp\{\omega\}\oplus
V_1\oplus\cdots\oplus V_{i+1})(\chi)$ which are not in $(\Fp\{\omega\}\oplus
V_1\oplus\cdots\oplus V_i)(\chi)$ is
$$
{pq^{i+1}-1\over p-1}-{pq^{i}-1\over p-1}\qquad
\left(\hbox{resp. }{q^{i+1}-1\over p-1}-{q^{i}-1\over p-1}\right)
$$
according as $\chi=\omega$ or $\chi\neq\omega$. The ``\thinspace
level\thinspace'' of such a line $D$ is $d(D)=pi+j$, where $j\in[1,p[$ is
determined by
$$
\bar v(\chi)\equiv\bar v(\omega)-(pi+j)\pmod{p-1},
$$
and therefore depends only on $\bar v(\chi)$ and
$i\pmod{p-1}$~\footnote{(*)}{(2011/07/09) After having been so scrupulous
  about the twist $\{\omega\}$ by $\omega$ hitherto, it was inadvertently
  omitted at this point in the earlier versions.  This omission merely had the
  effect of permuting the terms.}.  This follows from the fact that the
$k[G]$-module $l[n]$ is the direct sum of $k$-lines $k\{\xi\}$ for those
characters $\xi:G\to\F_p^\times$ for which $\bar v(\xi)\equiv n\pmod{p-1}$
(prop.~32).

Hence the contribution of all $G$-stable lines for given~$i$, $\chi$, and $j$
(with $\bar v(\chi)\equiv\bar v(\omega)-(i+j)$) is
$$
p\left({q^{i+1}-q^{i}\over p-1}\right)q^{-(pi+j)}
$$
irrespective of whether $\chi=\omega$ or $\chi\neq\omega$, as we saw above.
As there are $p-1$ characters having a given ``\thinspace
valuation\thinspace'', the contribution of all $G$-stable lines for a given
$i$ is
$$
p(q^{i+1}-q^{i})q^{-pi}Q,\qquad Q=\sum_{j=1}^{p-1}q^{-j}.
$$
Now all that remains to be done in the characteristic-$p$ case is to sum over
all $i\in\N$ and use the fact that
$$
\sum_i(q^{i+1}-q^i)q^{-pi}=(q-1)\sum_i q^{i-pi}=(q-1){q^{p-1}\over
  q^{p-1}-1}={1\over Q},
$$
proving the formula in this case.  In the characteristic-$0$ case, the sum
extends only over $i\in[0,e[$ to give the contribution of {\it peu
  ramifi{\'e}es\/} extensions
$$
pQ\sum_{i=0}^{e-1}(q^{i+1}-q^i)q^{-pi}=p(1-q^{(1-p)e}).
$$
But to this must be added the contribution of the {\it tr{\`e}s
  ramifi{\'e}es\/} extensions, which correspond to lines in the $1$-eigenspace
${\cal F}_{pe}(1)=(K^\times\!/K^{\times p})(1)$ which are not in the
hyperplane ${\cal F}_{pe-1}(1)=\bar U_1(1)$.  The dimension of the ambient
space is
$$
\dim_{\Fp}(\Fp\{\omega\}\oplus V_1\oplus\cdots\oplus V_e\oplus\Fp)(1)
=\cases{
ef+1&if $1\neq\omega$,\cr
ef+2&if $1=\omega$,\cr}
$$
giving the following contributions in the two cases $1\neq\omega$,
resp. $1=\omega$~: 
$$
p\left({pq^{e}-1\over p-1}-{q^{e}-1\over p-1}\right)q^{-pe},\qquad
\left({p^2q^{e}-1\over p-1}-{pq^{e}-1\over p-1}\right)q^{-pe}.
$$
But these two are the same and equal $pq^{(1-p)e}$.  The total comes to~$p$
as expected, proving the degree-$p$ mass formula in characteristic~$0$ as
well.  (Note that our method computes the contributions $p.(1-q^{(1-p)e})$
and $p.q^{(1-p)e}$ of  {\it peu ramifi{\'e}es\/} and {\it tr{\`e}s
  ramifi{\'e}es\/} extensions separately, not just their sum~$p$.)

This proof is reminiscent of the {\it
  F{\"u}hrerdiskriminantenproduktformel\/}, as applied before lemma~6 in
  \citer\further().  We thus obtain the following prime-degree case of
  Serre's mass formula~:

\th THEOREM 39 (Serre, 1978)
\enonce
Let\/ $F$ be a local field with finite residue field of characteristic~$p$
and cardinality~$q$.  When\/ $E$ runs through ramified separable degree-$p$
extensions of\/ $F$ (contained in a fixed separable algebraic closure
$\tilde F$ of\/~$F$), then
$$
\sum_E q^{-c(E)}
=p,
$$
where\/ $c(E)=v(\delta_{E|F})-(p-1)$, and $v(\delta_{E|F})>(p-1)$ is the
valuation of the discriminant\/ $\delta_{E|F}$ of\/ $E|F$ \footnote{\rm
  (*)}{\rm (2011/03/15) Clearly, our method can be used to count the number of
  degree-$p$ separable extensions $E$ of $F$ (contained in $\tilde F$) with a
  given $c(E)$, and also the number of $F$-conjugacy classes of such $E$.
  These numbers were determined by Krasner in all degrees~; using them, he
  gave a different proof of Serre's mass formula ({\it Remarques au sujet
    d'une note de J.-P.~Serre: ``\thinspace Une `formule de masse' pour les
    extensions totalement ramifi{\'e}es de degr{\'e} donn{\'e} d'un corps
    local\/\thinspace'' : une d{\'e}monstration de la formule de M.~Serre
    {\`a} partir de mon th{\'e}or{\`e}me sur le nombre des extensions
    s{\'e}parables d'un corps valu{\'e} localement compact, qui sont d'un
    degr{\'e} et d'une diff{\'e}rente donn{\'e}s}, Comptes Rendus {\bf 288}
  (1979)~18, pp.~A863--A865.)}.  \endth

The above proof is summarised in \S9.  Some refinements of the
degree-$p$ mass formula are also presented there.
 
\bigbreak

{\bf 8.  Tame extensions of prime degree}\pointir Let us end, for the sake of
completeness and contrast, with a word about the compositum $C'$ of all
degree-$p'$ extensions $E'$ of $F$, where $p'$ is a prime $\neq p$.  The
reader who has made it so far should have no difficulty in supplying proofs
modelled on \S4.  In any case, the assertions here are nothing but
translations of standard facts about tamely ramified extensions of~$F$ into
the language of~\S6~; they serve to further illustrate the general theory
of~\S4.

Let $K'=F(\zeta')$, where $\zeta'$ is a primitive $p'$-th root of~$1$,
$G'=\Gal(K'|F)$ and let $\omega':G'\to\F_{p'}^\times$ be the cyclotomic
character, so that $\sigma'(\zeta')=\zeta'^{\omega'(\sigma')}$ for every
$\sigma'\in G'$.  The extension $K'|F$ is unramified of degree equal to the
order of $q$ in $\F_{p'}^\times$, so $K'=F\Leftrightarrow
G'=\{1\}\Leftrightarrow\omega'=1\Leftrightarrow p'\,|\,(q-1)$.

The $\F_{p'}$-space $\overline{K'^\times}=K'^\times\!/K'^{\times p'}$ is of
dimension~$2$ and contains the line $\overline{{\goth o}'^\times}$ on which
$G'$ acts via $\omega'$, where ${\goth o}'$ is the ring of integers of $K'$~;
this line is canonically $G'$-isomorphic to $k'^\times\!/k'^{\times p'}$,
where $k'$ is the residue field of $K'$ and $G'$ has been identified with
$\Gal(k'|k)$.  The valuation provides a $G'$-isomorphism
$\overline{K'^\times}\!/\overline{{\goth o}'^\times}=\F_{p'}$, so that
$\overline{K'^\times}$ is isomorphic to $\overline{{\goth
    o}'^\times}\oplus\F_{p'}=\F_{p'}\{\omega'\}\oplus\F_{p'}$ as a filtered
$\F_{p'}[G']$-module.  Consequently, the only characters
$\chi':G'\to\F_{p'}^\times$ for which $\overline{K'^\times}(\chi')\neq\{1\}$
are $\chi'=\omega'$ and $\chi'=1$.

Every ramified $E'$ is of the form $F(\!\root p'\of\pi)$ for some uniformiser
$\pi$ of $F$, and hence $E'K'$ is cyclic (of degree~$p'$) over $K'$ (and
galoisian over $F$)~; it thereby gives rise to a $G'$-stable line
$D'\subset\overline{K'^\times}$.  Conversely, every $G'$-stable line
$D'\subset\overline{K'^\times}$ comes from some $E'$~; $\overline{{\goth
    o}'^\times}$ comes from the unramified $E'$.  If
$\chi':G'\to\F_{p'}^\times$ is the character through which $G'$ acts on $D'$,
then $D'$ comes from one (resp.~$p'$) $E'$ if $\chi'=\omega'$
(resp.~$\chi'\neq\omega'$).  More precisely, we have established a bijection
between the set of $G'$-stable lines in $\overline{K'^\times}$ and the set of
$F$-conjugacy classes of degree-$p'$ extensions of $F$.

Let $M'=K'(\!\root p'\of{K'^\times})$ be the maximal abelian extension of $K'$
of exponent $p'$.  It is easily seen that $C'=M'$.  Indeed, we have just seen
that $C'\subset M'$.  Next, $K'\subset C'$ because $K'$ is contained in the
galoisian closure of any ramified $E'$.  As $C'|K'$ is abelian of
exponent~$p'$, there is a subspace $T'\subset\overline{K'^\times}$ such that
$C'=K'(\!\root p'\of{T'})$.  As $T'$ contains every $G'$-stable line, we must
have $T'=\overline{K'^\times}$ and $C'=M'$.  

Life would be dry if everything had been so tame.
% It should get even wilder in degrees~$p^n$ ($n>1$).

{\it Remark 40}\pointir The foregoing can be used to prove Serre's mass
formula in degree~$p'$, just as \S6 was used in \S7.  Indeed, the only
character $\chi':G'\to\F_{p'}^\times$ for which there is a line
$D'\neq\overline{{\goth o}'^\times}$ in $\overline{K'^\times}(\chi')$ is
$\chi'=1$.  There are $p'$ (resp.~$1$) such $D'$ if $1=\omega'$
(resp.~$1\neq\omega'$) and each $D'$ comes from $1$ (resp.~$p'$) ramified
degree-$p'$ extension $E'$ of $F$ if $1=\omega'$ (resp.~$1\neq\omega'$).  This
balancing, similar to what we saw in \S7 for {\it tr{\`e}s ramifi{\'e}es\/}
extensions of degree~$p$, leads to the result.

\bigbreak

{\bf 9.  Some refined mass formulae}\pointir Let us summarise the
proof of the degree-$p$ mass formula and show how the same strategy
leads to certain refinements.  These consist mainly in computing the
contribution of those degree-$p$ separable extensions $E$ of $F$ for
which the group $G=\Gal(K|F)$ acts on the order-$p$ group $\Gal(EK|K)$
via a given character $G\to\F_p^\times$.

{\it Summary of the proof\/} (2011/07/03)\pointir There is a map
$E\mapsto EK$ (with $K=F(\!\root p-1\of{F^\times})$) from the set
${\cal S}_p(F)$ of degree-$p$ separable extensions of $F$ (in some
fixed separable algebraic closure of $F$) into the set ${\cal C}_p(K)$
of degree-$p$ cyclic extensions of $K$~; its image is the set ${\cal
C}_p(K,F)$ of those $L\in{\cal C}_p(K)$ which are galoisian over~$F$,
and if we go modulo the relation $\sim$ of $F$-conjugacy in ${\cal
S}_p(F)$, we get a bijection ${\cal S}_p(F)/\mathord\sim\to{\cal C}_p(K,F)$.

There is a natural bijection of ${\cal C}_p(K)$ with the set
$\P(K^\times\!/K^{\times p})$ or $\P(K^+\!/\wp(K^+))$ of lines in the $\F_p$-space
$K^\times\!/K^{\times p}$ or $K^+\!/\wp(K^+)$ respectively.  The subset ${\cal
  C}_p(K,F)$ corresponds to the set $\P(K^\times\!/K^{\times p})^G$ or
$\P(K^+\!/\wp(K^+))^G$ of $G$-stable lines for $G=\Gal(K|F)$.  Under the composite
map
$$
\Phi:{\cal S}_p(F)\to{\bf P}(\ )^G,\quad
\Phi(E)=D,\quad 
(EK=K(\root p\of D)\hbox{ or }EK=K(\wp^{-1}(D))),
$$
the set ${\cal C}_p(F)$ of degree-$p$ cyclic extensions of $F$ is in
bijection with the set ${\bf P}((\ )(\omega))$ of lines in the
$\omega$-eigenspace for the action of~$G$, where $\omega:G\to\F_p^\times$ is
the mod-$p$ cyclotomic character (resp.~$\omega=1$ is the trivial character,
so that $K^+\!/\wp(K^+)(\omega)=F/\wp(F)$).  On the complement of ${\bf P}((\ 
)(\omega))$, the fibres of the map $\Phi$ have $p$ elements each, mutually
conjugate as extensions of $F$.  Note finally that $\Phi$ sends the unramified
degree-$p$ extension $F_p$ of $F$ to the line $\bar U_{pe}$ in characteristic~$0$
(where $e$ is the absolute ramification index of $F$) or $\bar\ogoth$ in
characteristic~$p$.  Here, $\bar\ogoth$ is the image of $\ogoth$ in
$K^+\!/\wp(K^+)$, where $\ogoth$ is the ring of integers of $K$, and $\bar U_{pe}$
is the image of $U_{pe}=1+\pgoth^{pe}$ in $K^\times\!/K^{\times p}$, where
$\pgoth$ is the unique maximal ideal of $\ogoth$.

Some of these facts are summarised in the following commutative diagram
$$\def\rrr#1{\hfl{#1}{}{5mm}}\def\ddd#1{\vfl{#1}{}{5mm}}
\diagram{
F_p&\in&{\cal C}_p(F)&\rrr{\subset}&{\cal S}_p(F)
 &\rrr{\hbox{``$p:1$''}}&{\cal C}_p(K,F)&\rrr{\subset}&{\cal C}_p(K)\cr
\ddd{}&&\ddd{\sim}&&&&\ddd{\sim}&&\ddd{\sim}\cr
(\bar U_{pe}\hbox{ or }\bar\ogoth)&\in&{\bf P}((\
  )(\omega))&\rrr{}&\rrr{\subset}&\rrr{\hbox{\phantom{``$p:1$''}}} &{\bf P}(\ )^G&\rrr{\subset}&{\bf
  P}(\ )\cr }$$ in which $(\ )$ is to be replaced by
  $K^\times\!/K^{\times p}$ or $K^+\!/\wp(K^+)$ respectively, and ``$p:1$''
  means that the map is $p$-to-$1$ except on the subset ${\cal
  C}_p(F)$, on which it is bijective.  In other words, every $D\in{\bf
  P}(\ )^G$ has $p$ preimages in ${\cal S}_p(F)$ unless the character
  through which $G$ acts on $D$ is $\omega$, in which case there is
  only one preimage.

The measure $c(E)=v(\delta_{E|F})-(p-1)$ of wild ramification of a ramified
$E\in{\cal S}_p(F)$ equals $d(D)$, the ``\thinspace level\thinspace'' of the
line $D$ corresponding to $E$ --- the first step in the filtration on
$K^\times\!/K^{\times p}$ or $K^+\!/\wp(K^+)$ to which $D$ belongs, when the
``\thinspace level\thinspace'' of $\bar U_{pe}$ and $\bar\ogoth$ is declared
to be~$0$ and the increasing filtration is indexed by $[0,pe]$
(resp.~$\N$).  In other words, the stratification on ${\cal S}_p(F)$
given by $c$ corresponds to the stratification on ${\bf P}(\ )^G$
coming from the filtration on $K^\times\!/K^{\times p}$ or $K^+\!/\wp(K^+)$ respectively.

The filtered $\F_p[G]$-module $K^\times\!/K^{\times p}$ or $K^+\!/\wp(K^+)$
was shown in~\S6 to be isomorphic to
$$
\Fp\{\omega\}\oplus l[e-1]\oplus\cdots 
\oplus l[e-b^{((p-1)e)}]\oplus\Fp
$$
or to $\Fp\oplus l[-1]\oplus\cdots\oplus l[-b^{(n)}]\oplus\cdots$
respectively, where $l=\ogoth/\pgoth$ is the residue field of $K$, $l[i]$
denotes the $k[G]$-module $\pgoth^i/\pgoth^{i+1}$ (where $k$ is the residue
field of $F$) and $\{\omega\}$ denotes twist by the character~$\omega$.  Also,
$b^{(n)}$ is the $n$-th \hbox{prime-to-$p$} integer (so that
$b^{((p-1)e)}=pe-1$)~; here, $n$ runs over all integers $>0$ in
characteristic~$p$ but only from $1$ to $(p-1)e$ in characteristic~0, which is
also when the last term $\Fp$ appears, accounting for {\it tr{\`e}s
  ramifi{\'e}es\/} extensions.  

Put $V_{i+1}=l[-(ip+1)]\{\omega\}\oplus\cdots\oplus l[-(ip+p-1)]\{\omega\}$.
The filtered $\Fp[G]$-module $K^\times\!/K^{\times p}$ is isomorphic to
$\Fp\{\omega\}\oplus V_1\oplus\cdots\oplus V_e\oplus\Fp$ in characteristic~$0$
(resp.~$K^+\!/\wp(K^+)$ is isomorphic to $\Fp\oplus V_1\oplus V_2\oplus\cdots$ in
characteristic~$p$).

Everything is now in place for proving the degree-$p$ mass formula by
rewriting it as an appropriate sum over $D\neq\bar U_{pe}$
(resp.~$D\neq\bar\ogoth$) in ${\bf P}(K^\times\!/K^{\times p})^G$ (resp.~${\bf
  P}(K^+\!/\wp(K^+))^G$) as in \S7.  This is done by successively computing the
contribution of $G$-stable lines in $\Fp\{\omega\}\oplus V_1\oplus\cdots
\oplus V_{i+1}$ which are not in $\Fp\{\omega\}\oplus V_1\oplus\cdots\oplus
V_i$, and adding them all up for $i\in[0,e[$ in characteristic~$0$
(resp.~$i\in\N$ in characteristic~$p$). Finally, in characteristic~$0$, one
adds the contribution of $G$-stable lines in $K^\times\!/K^{\times p}$ which
are not in $\bar U_1$.

\smallbreak

{\it Remark 41\/} (2011/07/14)\pointir As the ``\thinspace level\thinspace''
of the line $\bar U_{pe}$ or $\bar\ogoth$ is~$0$, logic dictates that we pose
$c(F_p)=0$ for the {\it unramified\/} degree-$p$ extension $F_p$ of $F$.  The
degree-$p$ mass formula then becomes $\sum_{E\in{\cal S}_p(F)}q^{-c(E)}=1+p$.
More interestingly, there is a map ${\cal S}_p(F)\to F^\times\!/F^{\times
  p-1}$ which sends every $E$ to the character in
$\Hom(G,\F_p^\times)=F^\times\!/F^{\times p-1}$ through which $G$ acts on the
$G$-stable $\Fp$-line $D\in{\bf P}(\ )^G$ corresponding to $E$.  Our analysis
makes it possible to compute $\sum_{E\mapsto\chi}q^{-c(E)}$ for any given
$\chi\in F^\times\!/F^{\times p-1}$.

Fix a character $\chi:G\to\F_p^\times$ and an index $i\in[0,e[$
(resp.~$i\in\N$), and define $j_{\chi,i}\in[1,p[$ by the requirement $\bar
v(\chi)\equiv\bar v(\omega)-(i+j_{\chi,i})\pmod{p-1}$.  We have seen that
every $\Fp$-line
$$
D\subset(\Fp\{\omega\}\oplus V_1\oplus\cdots\oplus V_{i+1})(\chi),\quad
D\not\subset(\Fp\{\omega\}\oplus V_1\oplus\cdots\oplus V_i)(\chi),
$$
has ``\thinspace level\thinspace'' $pi+j_{\chi,i}$ (cf.~prop.~32) and their
contribution to the mass formula is
$$
p\left({q^{i+1}-q^{i}\over p-1}\right)q^{-(pi+j_{\chi,i})}.
$$
So the contribution of $\chi$ --- it turns out to depend only on
$\bar v(\chi)$ --- is 
$$
{p(q-1)\over p-1}\sum_i q^{i-(pi+j_{\chi,i})}\qquad 
\cases{i\in[0,e[&if $\car(F)=0$,\cr i\in\N&if $\car(F)=p$,\cr}
$$
except when $\chi=1$ is the trivial character and $\car(F)=0$, in which
case {\it tr{\`e}s ramifi{\'e}es\/} extensions contribute a further
$p/q^{(p-1)e}$.  This sum can be evaluated upon remarking that
$j_{\chi,i}=j_{\chi,i'}$ if $i\equiv i'\pmod{p-1}$.

{\it Remark 42} (2011/07/14)\pointir Reconciling these results with the
formulae (1)--(3) in \citer\final(prop.~14--16) for $\chi=\omega$ amounts to
checking the identity
$$
pi+j_{\omega,i}=(p-1)b^{(i+1)}
$$
for all $i\in\N$.  Upon writing $i=(p-1)n_i+r_i$ (with $r_i\in[0,p-1[$,
$n_i\in\N$), we have $j_{\omega,i}=p-1-r_i$ and $b^{(i+1)}=i+1+n_i$, and hence
the identity.

{\it Example 43} (2011/10/06)\pointir Take $F=\F_3\series{\pi}$ so
  that $p=q=3$, $e=+\infty$, and $F^\times\!/F^{\times
  2}=\{\bar1,-\bar1,\bar\pi,-\bar\pi\}$.  For $\chi\in
  F^\times\!/F^{\times 2}$ and $i\in\N$,
$$
  \bar v(\chi)=\bar0\ \ \ \Longrightarrow\ \ \  j_{\chi,i}
  =\cases{2&if $i\equiv0\pmod2$\phantom,\cr1&if $i\equiv1\pmod2$,\cr}
$$
so the contribution of each of the unramified characters (namely
$\bar1$ and $-\bar1$) is $9/20$.  Similarly,
$$
  \bar v(\chi)=\bar1\ \ \ \Longrightarrow\ \ \  j_{\chi,i}
  =\cases{1&if $i\equiv0\pmod2$\phantom,\cr2&if $i\equiv1\pmod2$,\cr}
$$
so the contribution of each of the ramified characters (namely $\bar\pi$ and $-\bar\pi$) is $21/20$. 

{\it Example 44} (2011/10/07) Take $p=5$ and $F=k\series{\pi}$.  For
any given $w\in\Z/4\Z$, a character of ``\thinspace
valuation\thinspace'' $w$ contributes $5(q-1)CA_w/4$, where
$C=\sum_{a\in\N}q^{-16a}=q^{16}/(q^{16}-1)$ and
$$\eqalign{
A_{\bar 0}&=q^{-4}+q^{-7}+q^{-10}+q^{-13},\cr
A_{\bar 3}&=q^{-1}+q^{-8}+q^{-11}+q^{-14},\cr
A_{\bar 2}&=q^{-2}+q^{-5}+q^{-12}+q^{-15},\cr
A_{\bar 1}&=q^{-3}+q^{-6}+q^{-09}+q^{-16}.\cr
}$$
(2011/10/08) Using {\it Remark\/}~41, it is easy to determine
the contribution of cyclic extensions to the degree-$p$ mass formula
for {\it odd\/}~$p$ (the contribution is $2$ for $p=2$).  Let's do it
more generally for any any character of ``\thinspace
valuation\thinspace'' $\bar0$ first in the easier case
$F=k\series{\pi}$.

\th PROPOSITION 45 
\enonce
Suppose that\/ $p\neq2$ and that\/ $F=k\series{\pi}$, and let $\chi\in
F^\times\!/F^{\times p-1}$ be a character of ``\thinspace
valuation\thinspace'' $\bar v(\chi)\equiv0\pmod{p-1}$.  The
contribution of $\chi$ to the degree-$p$ mass formula is
$$
\sum_{E\mapsto\chi}q^{-c(E)}=
{pq^{p-2}(q-1)(q^{(p-2)(p-1)}-1)\over(p-1)(q^{p-2}-1)(q^{(p-1)^2}-1)}
$$
(sauf erreur\/).  In particular, for\/ $\chi=1$, this is the contribution of cyclic extensions.
\endth
It is a matter of evaluating the sum
$\sum_{i\in\N}q^{i-(pi+j_{\chi,i})}$ (cf.~{\it Remark\/}~41), where $j_{\chi,i}\in[1,p-1]$ is subject to 
$i+j_{\chi,i}\equiv0\pmod{p-1}$.  Write $i=(p-1)n+r$ ($r\in[0,p-2]$,
$n\in\N$), and notice that $j_{\chi,i}=p-1-r$, so that
$$
\eqalign{
\sum_{i\in\N}q^{i-(pi+j_{\chi,i})}&=\sum_{r=0}^{p-2}\left(\sum_{n\in\N}q^{-(p-1)^2n-(p-2)r-(p-1)}\right)\cr
&=q^{-(p-1)}\sum_{r=0}^{p-2}\left(q^{-(p-2)r}\sum_{n\in\N}q^{-(p-1)^2n}\right)\cr
&=q^{-(p-1)}{q^{(p-1)^2}\over
   q^{(p-1)^2}-1}\sum_{r=0}^{p-2}q^{-(p-2)r}\cr
&={q^{p-2}(q^{(p-2)(p-1)}-1)\over(q^{p-2}-1)(q^{(p-1)^2}-1)}.\cr
}$$
Plugging this value into the formula in {\it Remark\/}~41 gives the result.

\th COROLLARY 46
\enonce
For $p$ odd and\/ $F=k\series{\pi}$, the contribution of all
ramified\/ $E\in{\cal S}_p(F)$ for which there is an
unramified extension\/ $F'|F$ (depending on $E$) such that\/ $EF'|F$
be galoisian is\/
$$
{pq^{p-2}(q-1)(q^{(p-2)(p-1)}-1)\over(q^{p-2}-1)(q^{(p-1)^2}-1)}.
$$
\endth
Let $\chi:G\to\F_p^\times$ be the character through which $G$ acts on
the line $D\in{\bf P}(K^+\!/\wp(K^+))^G$ corresponding to $E$.  It suffices
to show that $EF'$ is galoisian over $F$ for some unramified extension
$F'$ of\/ $F$ if and only if $\bar v(\chi)\equiv0\pmod{p-1}$.

Now, in the notation of lemma~19, the galoisian closure of $E$ over
$F$ is $EF_\chi$ (where $F_\chi=K^{\Ker(\chi^{-1})}$).  But $F_\chi|F$
is unramified if and only if $\chi$ is unramified
($\Leftrightarrow \chi(G_0)=\{1\}\Leftrightarrow\bar v(\chi)\equiv0$),
for $\Gal(F_\chi|F)=\Im(\chi^{-1})$, so
$\Gal(F_\chi|F)_0=\chi^{-1}(G_0)$.  As there are $p-1$ unramified
characters, the corollary follows from prop.~44.

(2011/10/11) More generally, let $\chi\in F^\times\!/F^{\times p-1}$
be any character and let $a\in[0,p-2]$ be such that $\bar
v(\chi)\equiv-a\pmod{p-1}$.  Notice that if $i-a=(p-1)n+r$ for some
$r\in[0,p-2]$ and some $n\in\N$, then $j_{\chi,i}=p-1-r$, as follows
from the defining requirements $j_{\chi,i}\in[1,p-1]$ and $\bar
v(\chi)\equiv-(i+j_{\chi,i})\pmod{p-1}$. This simple device helps in
evaluating the sum $\sum_{i\in\N}$ by rewriting it as
$\sum_{i=0}^{a-1}+\sum_{r=0}^{p-2}\sum_{n\in\N}$, where the first
sum is empty if $a=0$ as in prop.~45.  The case of ramified characters
($a>0$) is treated in the following proposition.

\th PROPOSITION 47
\enonce
Suppose that\/ $p\neq2$ and that\/ $F=k\series{\pi}$.   For a
character\/ $\chi\in F^\times\!/F^{\times p-1}$ of\/ $G=\Gal(K|F)$, define\/
$a\in[0,p-2]$ by\/ $\bar v(\chi)\equiv-a\pmod{p-1}$.  Then
the contribution\/ $\sum_{E\mapsto\chi}q^{-c(E)}$ of\/ $\chi$ to the
degree-$p$ mass formula is\/ ${p(q-1)/(p-1)}$ times
$$
{q^{p-2}(q^{(p-2)a}-1)\over q^{(p-1)a}(q^{p-2}-1)}
 +{q^{p-2}(q^{(p-2)(p-1)}-1)\over q^{(p-1)a}(q^{p-2}-1)(q^{(p-1)^2}-1)}.
$$
\endth
[Notice that for $a=0$ the first term is $0$ and we retrieve the expression in prop.~45.]
We have to compute $\sum_{i\in\N}q^{i-(pi+j_{\chi,i})}$.  The first
term in the above expression corresponds to $i\in[0,a[$, interval in
which $j_{\chi,i}=a-i$. The second term corresponds to
$i\in[a,+\infty[$~; to evaluated it, write $i-a=(p-1)n+r$ (with
$r\in[0,p-2]$, $n\in\N$), and note
that $j_{\chi,i}=p-1-r$, as explained above.  So
$$
\eqalign{
\sum_{i\in[a,+\infty[}q^{i-(pi+j_{\chi,i})}
 &=\sum_{r=0}^{p-2}\left(\sum_{n\in\N}q^{-(p-1)^2n-(p-2)r-(p-1)(a+1)}\right)\cr
&=q^{-(p-1)(a+1)}\sum_{r=0}^{p-2}\left(q^{-(p-2)r}\sum_{n\in\N}q^{-(p-1)^2n}\right)\cr
&=q^{-(p-1)(a+1)}{q^{(p-1)^2}\over
   q^{(p-1)^2}-1}\sum_{r=0}^{p-2}q^{-(p-2)r}\cr
&={q^{p-2}(q^{(p-2)(p-1)}-1)\over q^{(p-1)a}(q^{p-2}-1)(q^{(p-1)^2}-1)}.\cr
}$$

{\it Remark 48} (2011/10/14)\pointir Evaluating the sum
$$
\sum_{a=0}^{p-2}{{(q^{(p-2)a}-1)(q^{(p-1)^2}-1)} +(q^{(p-2)(p-1)}-1)\over q^{(p-1)a}}
$$
gives $(q^{p-2}-1)(q^{(p-1)^2}-1)/q^{p-2}(q-1)$, which serves as a
check on the computations by reproving the degree-$p$ mass formula.

% {\it Remark 47} (2011/10/14)\pointir The identity
% $$
% \sum_{a=0}^{p-2}{(q^{(p-2)a}-1)(q^{(p-1)^2}-1) +(q^{(p-2)(p-1)}-1)\over q^{(p-1)a}}
% ={(q^{p-2}-1)(q^{(p-1)^2}-1)\over q^{p-2}(q-1)}
% $$
% serves as a check on the computations by reproving the
% degree-$p$ mass formula.

{\it Example 49} (2011/07/14)\pointir Take $F=\Q_3$, so that $p=q=3$, $e=1$,
and the respective contributions of the four characters in $F^\times\!/F^{\times
  2}=\{\bar1,-\bar1,\bar3,-\bar3\}$ are
$1+{1\over3}$, $1$, $1\over3$, $1\over3$.

{\it Example 50} (2011/07/14)\pointir Take $p=5$, $e=5$, and $\chi\in
F^\times\!/F^{\times4}$ with $\chi\neq1$ but $\bar v(\chi)=\bar0$.  Then
$j_{\chi,i}=1,4,3,2,1$ respectively for $i\in[0,5[$, and the contribution of
$\chi$ is $5(q-1)(q^{-1}+q^{-8}+q^{-11}+q^{-14}+q^{-17})/4$.  Also, the
contribution of $\omega$ is $5(q-1)(q^{-4}+q^{-7}+q^{-10}+q^{-13}+q^{-20})/4$.

(2011/10/16) Let us now treat the characteristic-$0$ case, so that $F$
is a finite extension of $\Q_p$ ($p\neq2$).  This case is different
from the characteristic-$p$ case in three ways~: $e$ is finite, $\bar
v(\omega)\equiv e$ may very well be $\not\equiv0\pmod{p-1}$, and
finally, there are {\it tr{\`e}s ramifi{\'e}es\/} extensions.
Nevertheless, as we have seen in \S6, there is enough structural
similarity between the $\F_p[G]$-module $K^+\!/\wp(K^+)$ in
characteristic~$p$ and the $\F_p[G]$-module $K^\times\!/K^{\times p}$
in characteristic~$0$ for us to carry out the same analysis.  The only
essential difference is that $\sum_{i\in\N}$ gets replaced by
$\sum_{i\in[1,e[}$.

As an illustration, we will compute the contribution of cyclic
extensions to the degree-$p$ mass formula.  More generally, we have 

\th PROPOSITION 51
\enonce
Suppose that\/ $F|\Q_p$ is a finite extension ($p\neq2$).  For
every\/ $\chi\in F^\times\!/F^{\times p-1}$ of ``\thinspace
valuation\thinspace'' $\bar v(\chi)\equiv\bar v(\omega)$, the contribution\/
$\sum_{E\mapsto\chi}q^{-c(E)}$ of\/ $\chi$ to the degree-$p$ mass
formula is
$$
{p(q-1)\over p-1}(I+J)
$$
(where $I$ and $J$ are described below) except for\/ $\chi=1$, when\/
{\rm tr{\`e}s ramifi{\'e}es\/} extensions contribute a further\/
$p/q^{(p-1)e}$.  In particular, these are the contributions of cyclic
extensions in the cases $\omega\neq1$, $\omega=1$ respectively.
\endth
We have to compute $\sum_{i\in[0,e[}q^{i-(pi+j_{\chi,i})}$.  Write
$e-1=(p-1)N+R$ (with $R\in[1,p-2]$ and $N\in\N$), and split the sum
over $i\in[0,e[$ into $i\in[0,(p-1)N[$ and $i\in[(p-1)N,e[$.  Notice
that if $i=(p-1)n+r$ (with $r\in[0,p-2]$ and $n\in\N$), then
$j_{\chi,i}=p-1-r$.  We have
$$
\eqalign{
I=\sum_{i=0}^{(p-1)N-1}q^{i-(pi+j_{\chi,i})}
&=\sum_{r=0}^{p-2}\left(\sum_{n=0}^{N-1}q^{-(p-1)^2n-(p-2)r-(p-1)}\right)\cr
&=q^{-(p-1)}\sum_{r=0}^{p-2}\left(q^{-(p-2)r}\sum_{n=0}^{N-1}q^{-(p-1)^2n}\right)\cr
&=q^{-(p-1)}{q^{(p-1)^2}\over q^{(p-1)^2}-1}{q^{(p-1)^2N}-1\over q^{(p-1)^2N}}\sum_{r=0}^{p-2}q^{-(p-2)r}.
}
$$
Notice that $I=0$ if $N=0$.  The second term is
$$
\eqalign{
J=\sum_{i=(p-1)N}^{e-1}q^{i-(pi+j_{\chi,i})}
&=\sum_{r=0}^R q^{-(p-1)^2N-(p-2)r-(p-1)}\cr
&=q^{-(p-1)^2N-(p-1)}\sum_{r=0}^R q^{-(p-2)r}.\cr
}
$$
It is left to the reader to do a similar analysis for other $\chi\in
F^\times\!/F^{\times p-1}$.

\th COROLLARY 52
\enonce
With the same notation, the contribution of those ramified\/ $E\in{\cal
S}_p(F)$ for which the maximal tamely ramified extension\/ $\tilde E_1$
of\/ $F$ in the galoisian closure\/ $\tilde E$ of\/ $E|F$ is\/
unramified over\/ $F$ is\/ $p(q-1)(I+J)$.
\endth
We know that the maximal tamely ramified extension of $F$ in $\tilde
E$ is $F_\chi$ (cf.~lemma~19), where $F_\chi$ is the fixed field of
$\Ker(\omega\chi^{-1})$ and $\chi:G\to\F^\times$ is the character
through which $G$ acts on the line $D\in\P(K^\times\!/K^{\times p})^G$
corresponding to $E$.  For $F_\chi$ to be unramified over $F$, it is
necessary and sufficient that $\omega\chi^{-1}$ be unramified, which
is equivalent to $\bar v(\chi)=\bar v(\omega)$.  As there are $p-1$
such $\chi$, the result follows from prop.~51.

{\it Remark\/ 53}\pointir Let us return to the general case of a finite
extension $F$ of ${\bf Q}_p$ or of ${\bf F}_p\series{\pi}$, where the
prime $p$ is odd, and fix an extension $F'$ of $F$ in $K$.  One can
similarly compute the contribution of all $E\in{\cal S}_p(F)$ such
that $EF'$ is galoisian over $F$~; the case of cyclic extensions
corresponds to the choice $F'=F$ and the degree-$p$ mass formula is
the case $F'=K$.  We could also require that $F'$ be the maximal
tamely ramified extension $\tilde E_1$ of $F$ in the galoisian closure
$\tilde E$ of $E$ over $F$.  The contributions in question are the sum
of the contributions of all $\chi\in F^\times\!/F^{\times p-1}$ such
that $F_\chi\subset F'$ (resp.~$F_\chi=F'$), in the notation of
lemma~19.  The basic case occurs when $F'|F$ is cyclic of degree
dividing $p-1$ (cf.~lemma~12 or {\it Remark\/}~23)~; the number of
such $F'$ can be easily computed \citer\hasse(Kap.~16)
or \citer\course(Lecture~18).

{\it Remark\/ 54} (2011/10/18)\pointir Similarly, given a group
$\Gamma$, extension of a subgroup of $\F_p^\times$ of order $n$ by a
group of order~$p$ (cf.~{\it Remark\/}~23), one can compute the
contribution of all $E\in{\cal S}_p(F)$ such that the group
$\Gal(\tilde E|F)$ is isomorphic to $\Gamma$.  As $\F_p^\times$ has
$\varphi(p-1)$ subgroups, there are $\varphi(p-1)$ possibilities for
$\Gamma$.  

{\it Example\/ 55} (2011/10/20)\pointir For example, when $n=2$, so
that $\Gamma$ is the dihedral group of order $2p$, and
$F=k\series{\pi}$, there are three characters in $\chi\in
F^\times\!/F^{\times p-1}$ for which $\Im(\chi^{-1})$ has order~$2$,
namely $\chi=\bar g^m, \bar\pi^m, (\bar g\bar\pi)^m$, where $g$ is a
generator of $k^\times$ and $m=(p-1)/2$.  The contribution of those
$E\in{\cal S}_p(F)$ whose galoisian closure $\tilde E$ over $F$ has
group $\Gal(\tilde E|F)$ isomorphic to $\Gamma$ is the sum of the
contributions of these three $\chi$, the first of which has
``\thinspace valuation\thinspace'' $\equiv0$ and the other two $\equiv
m\pmod{p-1}$.

{\it Remark\/ 56} (2011/07/14)\pointir Notice finally that when
$e<+\infty$, we can similarly compute the {\it number\/} of
$E\in{\cal S}_p(F)$ mapping to any given $\chi\in F^\times\!/F^{\times
p-1}$ under the composite
$$\def\rrr#1{\hfl{#1}{}{5mm}}
\gamma:{\cal S}_p(F)\to
{\bf P}(K^\times\!/K^{\times p})^G\to
\Hom(G,\F_p^\times)\rrr{\sim}
F^\times\!/F^{\times p-1}
$$
where the second map sends a $G$-stable line $D$ to the character through
which $G$ acts on $D$~; the cases $\chi=1,\omega$ would need a separate
treatement.  The isomorphism of $K^\times\!/K^{\times p}$ with
$\Fp\{\omega\}\oplus V_1\oplus\cdots\oplus V_e\oplus\Fp$ as a filtered
$\Fp[G]$-module is all that is needed.  The same method works for $e=+\infty$
if we restrict to $E$ with bounded $c(E)$~; here the case $\chi=1$ would be
special.

{\it Remark\/ 57} (2011/10/18)\pointir For a given cyclic extension
$F'$ of $F$ in $K$ of degree dividing $p-1$, one can similarly compute
the {\it number\/} of $E\in{\cal S}_p(F)$ whose galoisian closure
$\tilde E$ is $EF'$, or for which the group $\Gal(\tilde E|F)$ is
isomorphic to a given group $\Gamma$ as in {\it Remark\/}~54.

{\it Remark\/ 58} (2011/10/19)\pointir There are two naturally defined
functions on the set ${\cal S}_p(F)$ of degree-$p$ separable
extensions of $F$.  The first one $c:{\cal S}_p(F)\to\N$ is a measure
of wild ramification and the second one $\gamma:{\cal S}_p(F)\to
F^\times\!/F^{\times p-1}$ doesn't seem to have a name in the
literature.  I believe that these are the only natural functions on
this set, and that subsets defined in terms of $c$ and $\gamma$ are
the only natural subsets of ${\cal S}_p(F)$~\footnote{(*)}{This is not
literally true.  Some quadratic extensions of $\Q_2$ can be embedded
into quartic cyclic extensions whereas others cannot~; this behaviour
cannot be detected by $c$.  More generally, if $F^\times$ has an
element of order~$p$ (in which case $\omega=1$) but no element of
order~$p^2$, then some degree-$p$ cyclic extensions of $F$ (so
$\gamma=\omega$) can be embedded in a degree-$p^2$ cyclic extensions
and some cannot~; I don't think this can be detected by $c$ alone.}.
Be that as it may, the methods of this paper allow us to compute the
mass or the cardinal of any such subset.

\medbreak

{\bf 10. Acknowledgements}\pointir Intellectual debt to the authors of
\citer\delcorso() --- who had the basic idea of attaching lines $D\subset
K^\times\!/K^{\times p}$ to degree-$p$ extensions $E$ of $F$ in
characteristic~$0$ --- is once again gratefully acknowledged.  I warmly thank
Robin Chapman and Jack Schmidt for the proof of lemma~2, Joseph Oesterl{\'e}
for the proof of lemma~3 and for his critical remarks, and Jung-Jo Lee for a
thorough reading of the whole text.  Lastly, I should not forget R R Simha who
asked me (2009/09/24) if there was anything more to be done in local
arithmetic. 

\bigbreak
\unvbox\bibbox 

\bye